\documentclass{article}
\usepackage{amsmath,amssymb,graphicx,stmaryrd}

\newtheorem{theorem}{Theorem}

\newtheorem{proposition}{Proposition}
\newtheorem{corollary}{Corollary}
\newtheorem{lemma}{Lemma}
\begin{document}
\title{ Generating Functionals of Random Packing Point Processes: From Mat\'ern to $\infty$-Mat\'ern }
\author{ Nguyen Tien Viet and Fran\c{c}ois Baccelli}
\maketitle
\begin{abstract}
 In this paper we study the generating functionals of several random packing processes: the classical Mat\'ern hard-core model; its extensions, the $k$-Mat\'ern models and the $\infty$-Mat\'ern model,  which is an example of random sequential packing process. We first give a sufficient condition for the $\infty$-Mat\'ern model to be well-defined (unlike the other two, the latter  may not be well-defined on unbounded spaces). Then the generating functional of the resulting point process is given for each of the three models as the solution of a differential equation. Series representations and bounds on the generating functional of the packing models are also derived. Last but not least, we obtain moment measures and Palm distributions of the considered packing models departing from their generating functionals.
\end{abstract}

\section{Introduction}\label{S:Intro}
Random packing models (RPM) are point processes (p.p.s) where points "contending" with each other cannot be simultaneously present. These p.p.s play an important role in many studies in physics, chemistry, material science, forestry and geology. The first use of RPMs is to describe systems with hard-core interactions. The most important applications are reactions on polymer chains \cite{GHH74}, chemisorption on a single-crystal surface \cite{E93}, and absorption in colloidial systems \cite{SVS00}. In these models, each point (molecule, particle,$\cdots$) in the system occupies some space, and two points with overlapping occupied space contend with each other. Another example is the study of seismic and forestry data patterns \cite{VEFSC07}, where RPMs are used as a reference model for the data sets under consideration.\\
Recently, the study of wireless communications by means of stochastic geometry gave rise to another type of application. In wireless communications, each point (node, user, transmitter,$\cdots$) does not occupy space but instead generates interference to other points \cite{HGA10,BB09}. Two points contend if either of them generates too much interference to the other. The present paper is mainly motivated by this kind of application. Nevertheless, the models we consider here are quite general and can be applied to other contexts. In particular, we will study here:
\begin{itemize}
	\item The Mat\'ern ``hard-core`` model: this is one of the classical models in the study of packing processes. Defined by Mat\'ern in \cite{M60} as the type II model and now usually  referred to as the Mat\'ern hard-core model, one can think of it as a dependent thinning of a Poisson p.p. $\Phi$, which is called the proposed p.p. The contention between two proposed points is determined by a \emph{contention mechanism}. In the original setting, this consists in letting each proposed point $x$ occupy a ball $B(x,r)$ of radius $r$ centered at $x$. Two points with overlapping balls, or equivalently, with Euclidian distance smaller than $2r$, contend with each other. Once the contention between points is determined, a \emph{retention mechanism} is used to prohibit the simultaneous presence of any two contending points. This mechanism assigns a uniform independent random mark $t(x)$ in $[0.1]$ to each proposed point $x$. A point is retained iff it has the smallest mark among its contenders.\\
	There are many extensions of the original Mat\'ern model, all of which consist in generalizing the contention mechanism whilst keeping the retention mechanism unchanged. The first extension of this kind is the "soft-core" model which randomizes the parameter $r$. Each proposed point $x$ occupies a ball $B(x,r(x))$ of random radius $r(x)$. Two points $x$ and $y$ contend iff the balls $B(x,r(x))$ and $B(y,r(y))$ overlap. A further extension in this direction leads to the contention mechanism based on the Boolean germ-grain model associated with $\Phi$, where two proposed points contend iff their grains overlap. In this paper, we consider an even more general setting: the \emph{random conection model} \cite{FBCMB05}. In this model, the contention between two points $x$ and $y$ is determined by a symmetrical random Boolean contention field $C(x,y)$: "$x$ and $y$ contend with each other iff $C(x,y)=1$". One can easily see that this model admits all the above instances of Mat\'ern models as special cases. For example, by letting $C(x,y)=\textbf{1}_{|x-y|<2r}$, we find back the original Mat\'ern model.\\
	There are many interpretations for the retention mechanism. One may interpret the marks $t(.)$ as the priority given to points, where the higher the priority, the smaller the mark. The retention mechanism retains the points with highest priority among its contender. One may also think of $t(.)$ as the time when a point is proposed to the p.p. This gives rise to the following temporal point of view, which will be used frequently in this paper. At anytime $t$ from $0$ to $1$, a fresh Poisson p.p. with some infinitesimal rate arrives. Upon arrival, a point $x$ is discarded if it finds any contender which has already arrived; otherwise it will be retained. There is also a graph theoretic construction, which is very useful in proving the existence and the relationships between the models. Given a proposed p.p. $\Phi$ and a contention field, we first construct the conflict graph of $\Phi$, which has the points of $\Phi$ as vertices. A directed edge is added from $x$ to $y$ if $x$ contends with $y$ and $t(x)>t(y)$. It is not difficult to see that the conflict graph is acyclic and that the Mat\'ern model retains all the points with out-degree $0$. We will come back to this construction in Subsection \ref{S:CS well defined}\\
    
	\item The $k$-Mat\'ern models: It has been observed that the Mat\'ern model is quite conservative. Roughly speaking, if the purpose of packing models is to retain from the process of proposed points a subset which is contention free, the subset selected by the Mat\'ern model is not maximal in the set theoretic sense.
	\begin{figure}
	   \centering
	     \includegraphics[width=14 cm,height=2.2 cm]{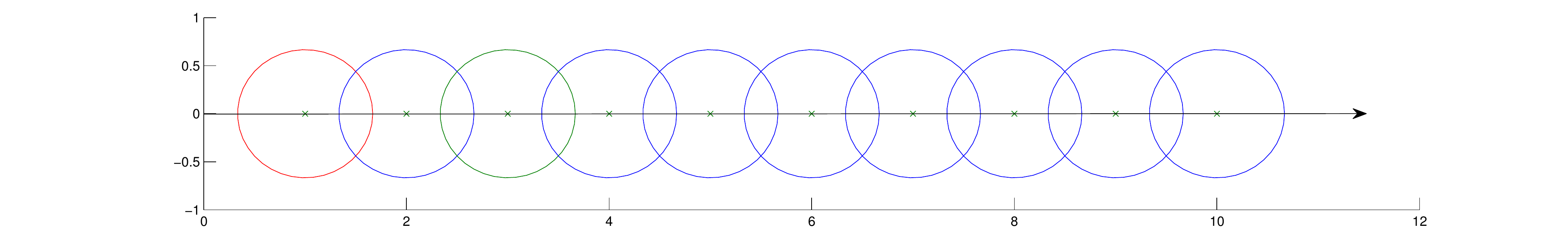}
% Pack toy example.eps: 0x0 pixel, 300dpi, 0.00x0.00 cm, bb=
	      \caption{Packing toy example. The points are at position $1,2,\cdots$, each of them occupies a circle of radius $2/3$. $t$ marks of $i$ is $1-i^{-1}$. The Mat\'ern model only retains $1$ while the $2$-Mat\'ern model retains $1$ and $3$.}
	    \label{fig:Packing toy example}
	  \end{figure}

	To see this more clearly, let us consider the following toy example (fig.\ref{fig:Packing toy example}). There is a sequence of points at the positions $1,2,\cdots$ on $[0,\infty)$. Each point occupies a ball of radius $2/3$ around it. The mark $t$ is assigned to the point at position $i$ is $1-1/i$. Thus, in this setting, the retention mechanism of the Mat\'ern model only retains the point at position $1$ whereas either all points at odd positions or all points at even positions form a contention free set.\\
	 In practice, A. Busson and G. Chelius \cite{BCG09} have conducted a simulation where the contention mechanism described in the above example is used and the Mat\'ern retention mechanism is applied to a Poisson point process. Let $r$ be the radius of the disc occupied by points in the Poisson p.p. Then, for each retained point, the disc of radius $2r$ centered on this point is its blocked area, i.e. no other point can be retained within this area. In the jamming regime, i.e. when the intensity of the proposed Poisson p.p. goes to $\infty$, only $75\%$ of the plane is blocked. As the intensity of the proposed p.p. is very large, one can find a proposed point in the un-blocked area with high probability. This point is not retained by Mat\'ern model although it does not contend with any retained point.   \\
	To alleviate this drawback, we propose the $k$-Mat\'ern models. Let us start with the example of fig.\ref{fig:Packing toy example}. The point at position $3$ is prohibited  by (i.e. contends with and has smaller $t$ mark than) the point at position $2$, which has already been prohibited by a retained point. Thus, retaining the point at position $3$ can do no harm. Generalizing this, we can retain more points by allowing any point that is not prohibited by any point retained by the Mat\'ern model. This results in the $2$-Mat\'ern model. However, this model also has its own drawback; it is ''over-greedy`` since two points retained by the $2$-Mat\'ern model can still contend with each other. For example, consider the situation in fig. \ref{fig:greedyexample}. There are five points at five vertices of a pentagon. Any two adjacent vertices contend with each other. The $t(.)$ marks are: $t(A)=.1$, $t(B)=.2$, $t(C)=.4$, $t(D)=.5$ and $t(E)=.3$. Thus, the Mat\'ern model retains vertex $A$, and the $2$-Mat\'ern model retains $A,D$ and $C$ (because D and C do not contend with $A$, which is the only point retained by the Mat\'ern model). However, we can see that $D$ and $C$ contend since they are adjacent. \\
\begin{figure}
 \centering
 \includegraphics[width=10cm]{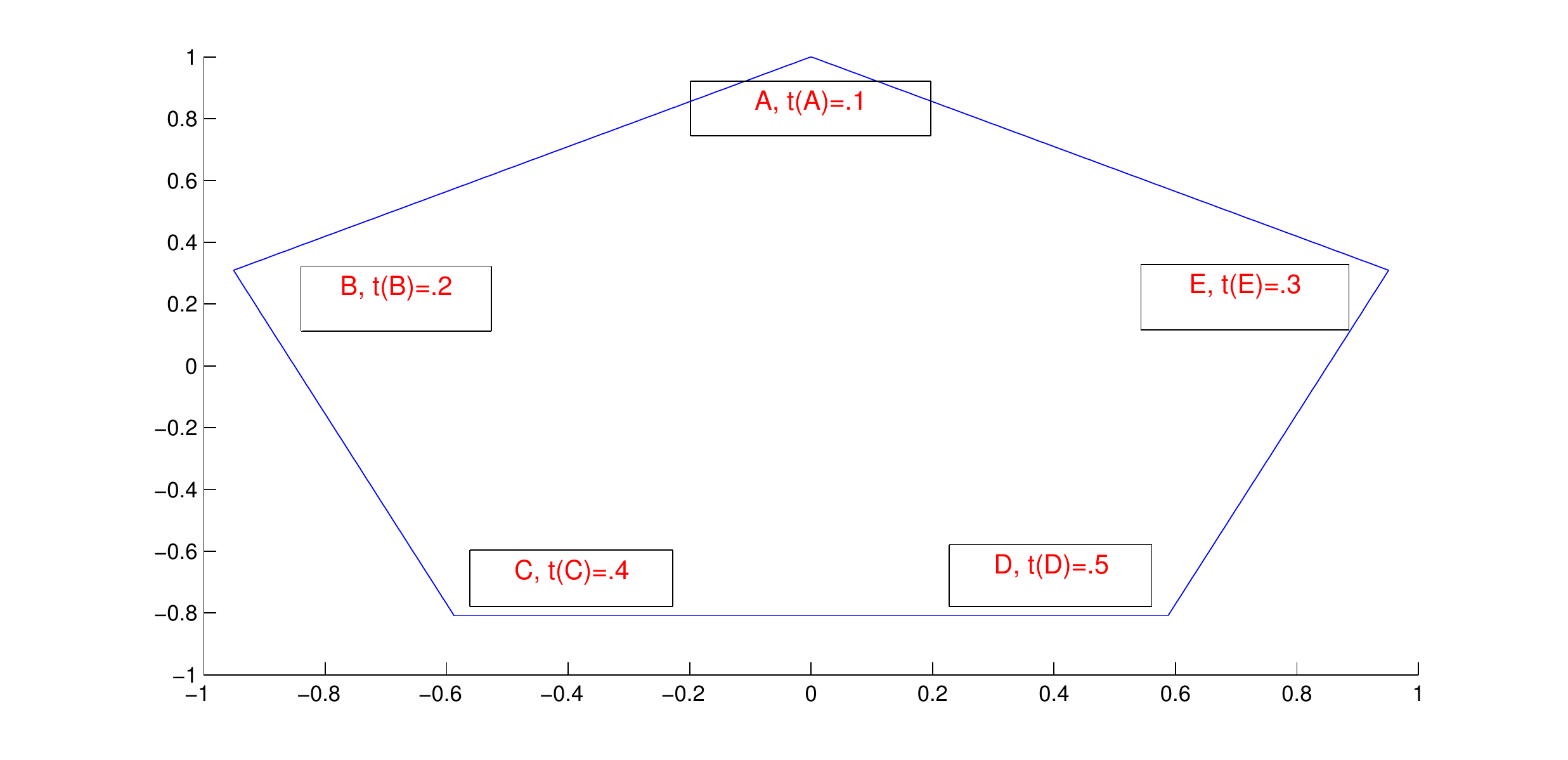}
 % greedyexample.eps: 0x0 pixel, 300dpi, 0.00x0.00 cm, bb=  -58   248   654   593
 \caption{An example where two points retained by the $2$-Mat\'ern model contend. Five points are placed at the vertices of a pentagon. Two adjacent vertices contend with each other. The $t$ marks are as shown in the figure.}
 \label{fig:greedyexample}
\end{figure}

	By iterating the procedure used to obtain the $2$-Mat\'ern model, we can recursively define the $k+1$-Mat\'ern model from the $k$-Mat\'ern model as follows: the subset of points retained by the $k+1$-Mat\'ern model is the proposed points which are not prohibited by any point retained by the $k$-Mat\'ern model. Within this context, the Mat\'ern model is the $1$-Mat\'ern model and the proposed p.p. $\Phi$ is the $0$-Mat\'ern model.\\
	More formally, the $k$-Mat\'ern model is defined as a dependent thinning of the proposed Poisson p.p. $\Phi$. The contention mechanism is the same as that of the Mat\'ern model. As for the retention mechanism, each proposed point $x$ is assigned a uniform independent random mark $t(x)$ in $[0,1]$ and a sequence of binary marks $e_{0}(x),e_{1}(x),\cdots$. The meaning of this sequence is: ''$x$ is retained by the $k$-Mat\'ern model iff the mark $e_{k}(x)=1$''. In the temporal view-point, at each time $t$ from $0$ to $1$, a fresh Poisson p.p. arrives. Upon arrival, each point $x$ is retained by the $k+1$-Mat\'ern model ($e_{k+1}(x)=1$) iff it finds no arrived contender which is retained by the $k$-Mat\'ern model. In the graph theoretic construction, we construct the conflict graph as for the Mat\'ern model. A point $x$ is assigned the mark $e_{k+1}(x)=1$ if for every child $y$ of $x$ (vertex that has an edge from $x$ to it, we recall that the conflict graph is acyclic so the notion of child is well-defined), we have $e_{k}(y)=0$. Note that $e_{0}(x)=1$ for all $x$ in $\Phi$ by definition.\\
	Before going to our next packing model, notice that the even-Mat\'ern models are not, strictly speaking, ``packing'' models. As we have already seen with the $2$-Mat\'ern model, two points retained by an even-Mat\'ern model can still contend with each other, so that the retained subset is not necessarily contention free. In general, the odd-Mat\'ern models retain points in a conservative way while the even-Mat\'ern ones do so in an over-greedy way.\\

	\item $\infty$-Mat\'ern model: The last model considered in this paper is the one arising in the study of wireless ad hoc networks using the Carrier Sensing Multiple Access (CSMA) protocol \cite{BB09}. This model is named $\infty$-Mat\'ern model, as we will show in Section \ref{S:CS well defined} that it is the ``limit'' of the $k$-Mat\'ern models when $k$ goes to $\infty$. It is also a ``perfect'' packing model as its retention mechanism is neither too conservative nor too greedy.  \\
	
	The $\infty$-Mat\'ern model can also be seen as a dependent thinning of the proposed Poisson p.p. $\Phi$. This thinning shares the same contention mechanism as all $k$-Mat\'ern models. The only difference lies in the retention mechanism. At this point, it would be convenient to keep the notation analogous to that of $k$-Mat\'ern models. Each proposed point $x$ is assigned a uniform mark $t(x)$ on $[0,1]$ and a Boolean mark $e_{\infty}(x)$ whose meaning is: ``$x$ is retained by the $\infty$-Mat\'ern models iff $e_{\infty}(x)=1$``. Temporally, when the point $x$ arrives at time $t(x)$, it is assigned the mark $e_{\infty}(x)=1$ (retained) iff it finds no contender which has arrived before and has already been assigned mark $1$. If one considers only the p.p. retained by the $\infty$-Mat\'ern model, it falls into a class of processes called Spatial Pure Birth p.p.s, which was suggested by M\o{}ller ( see \cite{SS00}, Section 2, paragraph 2). A pure birth process in the plane is a Markovian system of p.p.s indexed by time $t$, or equivalently, a Markov process whose states are p.p.s. At an instant of time, points are "born", i.e. added to the existing p.p. where the birth probability only depends on the current configuration $\phi$ of the process. The process is controlled by the birth rate $b$, which is a positive function $b(x,\phi)$ satisfying:
	\begin{align*}
	\int_{A}b(x,\phi)dx <\infty,
	\end{align*}
	for any bounded subset $A$ of $\mathbb{R}^{2}$ and any configuration $\phi$. Thus, given that the configuration of the process at time $t$ is $\phi$, the probability that there is a point born in $A$ in the infinitesimal time interval $[t,t+ds]$ is: 
	\begin{align*}
	ds\int_{A}b(x,\phi)dx+o(ds).
	\end{align*}
	In the $\infty$-Mat\'ern model, assuming that the intensity of the proposed process is $\lambda$, the birthrate can be expressed as:
	 \begin{align*}
	b(x,\phi) = \lambda \textbf{P}(\textnormal{a point at $x$ does not contend with any point in $\phi$}).
	\end{align*} 
	The graph-theoretic definition of the retention mechanism of the $\infty$-Mat\'ern is as follows: a point is retained iff, in the conflict graph, none of its children is retained. This is a recursive definition, the well-definedness of which depends on the structure of the conflict graph. In Section \ref{S:CS well defined}, we will show that under some mild conditions, the conflict graph of $\Phi$ is such that the previous recursive definition works well.    
\end{itemize}
\section{State of the Art}
Due to their wide range of applications, random packing models have attracted a lot of attention. The methods employed range from real experiments \cite{Feder80,ASZB94}, to simulations \cite{BS70,BCG09} and numerical approximations \cite{CaserHilhorst95,BCG09}.\\
On the other hand, the vast number of experimental publications here is in sharp contrast with the lack of rigorous mathematical results, especially for dimensions larger than 2. In dimension 1, packing models are more analytically tractable due to the shielding effect, \cite{SW86}. The first model of this kind is the car parking model which was independently studied by A. Renyi \cite{R58} and by A. Robin and H. Dvoretzky \cite{DR64}. In this model, cars of fixed length are parked in the same manner as in the $\infty$-Mat\'ern model. Consider an observation window $[0,x]$ and let $N(x)$ be the number of cars parked in this window when there is an infinite number of cars to be parked (saturated regime). R\'enyi showed that $N(x)$ satisfies the law of large number (LLN):
\begin{align*}
 \lim_{x\rightarrow \infty} \frac{N(x)}{x}= C \approx 0.74759 \qquad \textnormal{a.s.},
\end{align*}
where $C$ is called the packing density.
 A. Robin and H. Dvoretzk \cite{DR64} sharpened this result to a central limit theorem (CLT):
\begin{align*}
 \frac{N(x)-Cx}{\sqrt{\rm{var}(N(x))}}\longrightarrow\mathcal{N}(0,1) \qquad \textnormal{ in distribution as } x\rightarrow \infty.
\end{align*}
 Various extensions of the above models were considered like the non-saturated regime (the number of cars to be parked is finite), random car lengths, etc. \cite{Mu68,CFJP98}. The latter is also known under the name random interval packing  and has many applications in resource allocation in communication theory. For the above models, the obtained results concern the packing density, LLN, CLT, the distribution of packed intervals and that of vacant intervals.\\
	For dimension more than $1$, in his PhD thesis \cite{M60}, B. Mat\'ern introduced 3 types of hard-core models. Among these, the Mat\'ern type I model is the simplest. It is constructed as a dependent thinning of a Poisson proposed p.p. Each proposed point has a hard-disk of radius $r$ attached to it. A point is accepted if its attached disk does not overlap with that of any other point. The Mat\'ern type II models have already been discussed in Section 1 and can be seen as an instance of the 1-Mat\'ern model. The Mat\'ern type III model, which Mat\'ern just briefly mentioned, can be seen as a special case of the $\infty$-Mat\'ern model in this paper. Moment measures and correlation functions can be computed for the Mat\'ern type II models. The Mat\'ern type III model proves to be far less tractable. In his work, Mat\'ern asserted:''even an attempt to find the [packing density] tends to rather formidable mathematics``.\\
 At the same time, I. Palasti \cite{P60} considered an extension of the car parking problem to the plane, where  cars are rectangles of fixed sides. The car parking problem for the $n$ dimensional space is then defined in the same manner. She conjectured that the packing density is $C^{n}$ for the car parking model in $n$ dimensional space. This conjecture was, however, invalidated later by experimental data \cite{BS82}.\\
 One may recognize that the Mat\'ern type III model and the car parking model are just two names of the same object. In fact, this model is most well-known as Random Sequential Absorbtion model (RSA). If we consider the $\infty$-Mat\'ern model defined only in a finite window, then it coincides with the Poisson RSA model. The reader should refer to the paper by H. Solomon and H. Weiner \cite{SW86} for a review on the progress in the study of the RSA model up to that time.\\
The most noticable advance in this field is a series of papers by M.D. Penrose, J.E. Yukich and Y. Baryshnikov. Based on a general LLN and CLT for \emph{stablizing} functionals, the LLN and CLT were established for the $n$ dimensional RSA model in the non-saturated \cite{PMY02} regime. Y. Baryshnikov and J. E. Yukich \cite{BY02} later strengthened the above results by proving that, in the thermodynamic limit, the spatial distribution of the p.p. induced by the RSA model converges to that of a Gaussian field after a suitable rescaling in the non-saturated regime. The LLN and CLT for the $n$ dimensional RSA model under saturated regime is proven by T. Schreiber \emph{et al} \cite{SMPY07}.\\
It is worth noting here the difference between our approach and that of the works listed above. As expained earlier in Section \ref{S:Intro}, in the $\infty$-Mat\'ern model (which also admits the RSA model as a special case), there is  a time dimension, along which new points arrive to be packed. We are interested in the dynamic of the induced p.p.s as time passes and more points arrive. In contrast, time is fixed in the aforementioned works. The above authors consider the distribution of the points packed in an observed window and how it converges when the observed window grows to the whole space.  
\section{The Packing Models}\label{S:model}
Throughout this paper, every p.p. is equipped with independent and identically distributed (i.i.d.) marks which we call \emph{timers}. These timers are uniformly distributed on $[0,1]$. The timer of point $x$ is denoted by $t(x)$. \\
Let us begin with a preliminary transformation that will be used frequently.
For each p.p. $\Phi$, we define:
\begin{align}
  T_{A}(\Phi)= \left\{x \in \Phi\ \textnormal{ s.t. }\ t(x) \in A\right\},
\end{align}
for all measurable subsets $A$ of $\mathbb{R}^{+}$. This is the \emph{timer based thinning} 
of $\Phi$ in $A$.\\
The contentions of the points are encoded in the Boolean random field $C$.
The meaning of $C$ is that two points $x, y$ contend iff $C(x,y)=1$. Let $\textbf{P}(C(x,y)=1) =h(x,y)$.
The assumptions on $C$ are:
\begin{align*}
&  C(x,y)  = C(y,x)\ \forall x,y,\ \textnormal{(symmetry)};\\
&  C(x,y)  \textnormal{ are independent r.v.s for unordered pairs }(x,y);\\
&  h(x,y)=h(0,y-x)\textnormal{ } \forall x,y \textnormal{ (translation invariance)}. 
\end{align*}
Note that the contention field is defined independently of the timers. If $C(x,y)=1$ then $x$ and $y$ contend with each other regardless of their timers.
\subsection{Packing Models by Thinning}
In this subsection, we define the Mat\'ern model, the $k$-Mat\'ern model and the $\infty$-Mat\'ern model
as thinning transformations of another p.p.. These thinning transformations can be obtained by
iteratively applying the following transformation $\Psi$:
\begin{align*}
\Psi(\Phi_{1},\Phi_{2},C)= \{x \in \Phi_{1} \textnormal{ s.t. } \prod_{y \in T_{[0,t(x))}(\Phi_{2})}(1-C(x,y))=1 \}.
\end{align*} 
\begin{itemize}
	\item Mat\'ern model: the Mat\'ern model $M^{1}(\Phi,C)$ of $\Phi$ and of the contention field $C$, is obtained from the thinning transformation
	\begin{align}
 M^{1}(\Phi,C):=\Psi(\Phi,\Phi,C).\label{E:1 def}
	\end{align}
That is:
	\begin{align}
	 M^{1}(\Phi,C) = \{x \in \Phi \textnormal{ s.t. } e_{1}(x)=1\},
	\end{align}
where the thinning indicators $e_{1}(x)$ are defined by:
	\begin{align}
	e_{1}(x) = \prod_{y \in T_{[0,t(x))}(\Phi)}(1-C(x,y)) \label{E:1 thin ind}. 
	\end{align}

\item $k$-Mat\'ern model: the $k$-Mat\'ern models $M^{k}(\Phi,C)$ of $\Phi$ and $C$
 can be recursively defined as
	\begin{align}
	M^{k}(\Phi,C)= \Psi(\Phi,M^{k-1}(\Phi,C),C). \label{E:k def}
	\end{align}
One can also define the thinning indicator of the $k$-Mat\'ern model and write:
\begin{eqnarray}
M^{k}(\Phi,C) &= &\{x \in \Phi \textnormal{ s.t. } e_{k}(x)=1\}, \label{E:k redef}
\end{eqnarray}
with
\begin{eqnarray}
e_{k}(x) & = & \prod_{y \in T_{[0,t(x))}(M^{k-1}(\Phi,C))}(1-C(x,y))\nonumber\\
& = &\prod_{y \in T_{[0,t(x))}(\Phi)}(1-e_{k-1}(y)C(x,y)). \label{E:k thin ind}
\end{eqnarray}
	\item $\infty$-Mat\'ern model: the $\infty$-Mat\'ern model $M^{\infty}(\Phi,C)$ of $\Phi$ and $C$ can be defined as:
\begin{align}
	M^{\infty}(\Phi,C)= \Psi(\Phi,M^{\infty}(\Phi,C),C), \label{E:CSMA def}
	\end{align}
in a sense to be made clear in the next subsection. An equivalent expression is:
\begin{eqnarray}
M^{\infty}(\Phi,C)= \{x \in \Phi \textnormal{ s.t. } e_{\infty}(x)=1\}, \label{E:CSMA redef}
\end{eqnarray}
with $e_{\infty}$ satisfying the fixed point equation
\begin{eqnarray}
	e_{\infty}(x)& = & \prod_{y \in T_{[0,t(x))}(M^{\infty}(\Phi,C))}(1-C(x,y))
\nonumber\\
& = &\prod_{y \in T_{[0,t(x))}(\Phi)}(1-e_{\infty}(y)C(x,y)). \label{E:CSMA thin ind}
\end{eqnarray}
\end{itemize} 
In the above definitions, we adopt the convention that the product over an empty set is $1$. 
\subsection{The $\infty$-Mat\'ern Model is Well Defined}\label{S:CS well defined}
In this subsection, we give a sufficient condition for the $\infty$-Mat\'ern model to be well-defined.
We then prove that the $\infty$-Mat\'ern model is the limit of the $k$-Mat\'ern models in the sense that:
	\[
  \forall x\in \Phi \qquad e_{\infty}(x)=\lim_{k\rightarrow\infty}e_{k}(x) \textnormal{ a.s. }
  \]
  \subsubsection{Conflict Graph}
  Here we formally define the notion of conflict graph which will play a central role in the forthcoming proofs.
Given a p.p. $\Phi$ and a contention field $C$, the associated \emph{conflict graph} is the directed graph
$G(\Phi,C)=(\Phi,\mathcal{E}(\Phi))$, where $\mathcal{E}(\Phi)= \{(x,y) \in \Phi^{2}$  s.t. $t(y)<t(x), C(x,y)=1\}$.
That is, we put an edge from $x$ to $y$ iff $x$ and $y$ contend and the timer of $x$ is larger than that of $y$. It is easy to check that $G(\Phi,C)$ is acyclic.\\
We now introduce the notion of \emph{conflict indicator function} of a graph.
To this end we consider an arbitrary acyclic directed graph $H$ and a vertex $o$ in it.
Let $K^{*}_{H}(o)$ be the \emph{positive connected component} of $o$, i.e. the set of points $x$ 
such that there is a directed path in $H$ from $o$ to $x$. For all positive integer $n$, let $K^{n}_{H}(o)$ 
be the $n$-step \emph{positive connected component} of $o$, i.e. the set of points $x$ such that there is
a path of length at most $n$ from $o$ to $x$.  Note that $o$ is always contained in $K^{*}_{H}(o)$ and $K^{n}_{H}(o)$ by convention.
We say that $H$ has no infinite \emph{positive percolation} if we have $|K^{*}_{H}(x)|<\infty$ for all $x$ in $\Phi$.

In the following, by abuse of notation, we use the same notation for a set of vertices and the subgraph induced by this set. The conflict indicator function $e(H,o)$, which takes value in $\{0,1\}$, is recursively defined as:
  \begin{align}
  e(H,o):=\left\{
  \begin{array}{l r }
  1 &\textnormal{ if } H=(\{o\},\emptyset)\\
  \prod_{(o,x)\in H} (1-e(K^{*}_{H}(x),x)) &\textnormal{otherwise.} \label{E:Thin ind func}
  \end{array}
  \right.
  \end{align} 
 One can easily see that when $|K^{*}_{H}(o)|<\infty$, $e(H,o)$ is well-defined and equal to $e(K^{*}_{H}(o),o)$.
We will need the following lemmas in which we assume that the condition 
\begin{align}
 \int_{\mathbb{R}^{2}}h(0,x)dx=N <\infty \label{E:CSMA well def cond}
  \end{align}
 is satisfied:

\begin{lemma}\label{L:small finite per cond}
For an homogeneous Poisson p.p. $\Phi$ of intensity $\lambda<N^{-1}$, the conflict graph $G(\Phi,C)$
has no infinite \emph{positive} percolation a.s. 
\end{lemma}
\emph{Proof.}\\
For this proof we introduce the notion of infinite percolation. In $G(\Phi,C)$, we remove the direction of the edges and let $S(x)$ be the connected component of the resulting graph containing $x$. We say that the graph $G(\Phi,C)$ has no infinite percolation if $|S(x)|<\infty$ for all $x$. Clearly, if $G(\Phi,C)$ has no infinite percolation then it has no infinite positive percolation. We now show that $G(\Phi,C)$ has no infinite percolation a.s. For each $x$:
\begin{align*}
  &\textbf{E}[|S(x)|] \leq 1+ \sum_{n=1}^{\infty}\textbf{E}\left[\sum_{y \in \Phi}\textbf{1}_{\textnormal{there is a path of length $n$ from $x$ to $y$}}\right]\\
				&\leq 1+ \sum_{n=1}^{\infty}\textbf{E}\left[\textnormal{number of non self-intersecting paths of length $n$ starting from $x$}\right]\\
				&= 1+ \sum_{n=1}^{\infty}\textbf{E}\left[\sum_{y_{1},y_{2},\cdots,y_{n} \textnormal{pairwise different} \in \Phi}\textbf{1}_{C(x,y_{1})=1,C(y_{1},y_{2})=1,\cdots,C(y_{n-1},y_{n})=1}\right].
\end{align*}
Using the independence of $C(x,y_{1}),C(y_{1},y_{2}),\cdots,C(y_{n-1},y_{n})$ and the reduced Campbell formula we have:
\begin{align*}
  &\sum_{i=1}^{n}\textbf{E}\left[\sum_{y_{1},y_{2},\cdots,y_{n} \textnormal{pairwise different} \in \Phi}\textbf{1}_{C(x,y_{1})=1,C(y_{1},y_{2})=1,\cdots,C(y_{n-1},y_{n})=1}\right]\\
  &=\sum_{i=1}^{n}\int_{(\mathbb{R}^{2})^{n}}h(x,y_{1})h(y_{1},y_{2})\cdots h(y_{n-1},y_{n})\Lambda^{(n!)}(dy_{1},dy_{2},..dy_{n}),
\end{align*}
where $\Lambda^{(n!)}$ is the $n^{th}$ factorial moment measure of $\Phi$. This measure is known to be $\lambda^{n}l$, where $l$ denotes the Lebesgue measure on $\mathbb{R}^{2}$. Thus, using the fact that $h$ is translation invariant, we can change the variables from $y_{i}$ to $y_{i}-y_{i-1}$ and
$y^{1}$ to $y^{1}-x$ to get:
\begin{align*}
 &\textbf{E}[|S(x)|] \\
&\leq  1+\sum_{n=1}^{\infty}\int_{(\mathbb{R}^{2})^{n}}h(0,y_{1}-x)h(0,y_{2}-y_{1})\cdots h(0,y_{n}-y_{n-1})\lambda^{n}dy_{1}dy_{2}\cdots dy_{n}\\
		&= 1+\sum_{n=1}^{\infty}\int_{(\mathbb{R}^{2})^{n}}h(0,y_{1})h(0,y_{2})\cdots h(0,y_{n})\lambda^{n}dy^{1}dy^{2}\cdots dy^{n}\\
	       &= 1+\sum_{n=1}^{\infty}\lambda^{n}N^{n}
	       = (1-\lambda N)^{-1}
	       < \infty.
\end{align*}
Thus $|S(x)|<\infty$ a.s. and this proves the lemma.
\begin{flushright}
$\Box$
\end{flushright}

\begin{lemma}\label{L:Palm finite per cond}
If $G(\Phi,C)$ has no infinite positive percolation, then for all  $x \in \mathbb{R}^{2}$,
$G(\Phi\cup\{x\}^{+},C)$ has no infinite positive percolation either.
The $+$ superscript means that the timer for the additional point is set to $1$.
\end{lemma}
\emph{Proof.}\\
 Note that since (\ref{E:CSMA well def cond}) is satisfied, any point $x$ in the plane has only a finite number of contenders in $\Phi$ a.s.. Thus, by adding a point to $\Phi$ and letting the timer of this point be $1$, we can only add a finite number of edges to the conflict graph. These edges are the edges connecting the added point to its contenders in $\Phi$.
So, if $G(\Phi,C)$ has no infinite positive percolation,
then $G(\Phi\cup\{x\}^{+},C)$ has no infinite positive percolation neither.
\begin{flushright}
 $\Box$
\end{flushright}

\begin{lemma}\label{L:large finite per cond}
For an homogeneous Poisson p.p. $\Phi$ of intensity $\lambda<\infty$, the conflict graph $G(\Phi,C)$
has no infinite positive percolation.
  \end{lemma}
\emph{Proof.}\\
We prove this lemma by induction on the  intensity of $\Phi$.
Fix a positive $\lambda_{0}<N^{-1}$. If $\lambda<\lambda_{0}$,
then $G(\Phi,C)$ does not have infinite positive percolation thanks to
Lemma \ref{L:small finite per cond}.
Now let $n$ be a positive integer and suppose that $G(\Phi,C)$ does not have infinite positive percolation for any homogeneous Poisson p.p.
$\Phi$ with intensity smaller than $n\lambda_{0}$. Consider an homogeneous Poisson p.p.
$\Phi$ with intensity $\lambda$ such that $n\lambda_{0}\leq\lambda<(n+1)\lambda_{0}$.
Since $T_{[\frac{n\lambda_{0}}{\lambda},1]}(\Phi)$ is an homogeneous
Poisson p.p. with intensity $\lambda-n\lambda_{0}\le \lambda_{0}$, $G(T_{[\frac{n\lambda_{0}}{\lambda},1]}(\Phi),C)$
has no infinite percolation by Lemma  \ref{L:small finite per cond}.\\
Consider a point $x$ in $\Phi$. If $t(x)< \frac{n\lambda_{0}}{\lambda}$,
then 
$K^{*}_{G(\Phi,C)}(x)= K^{*}_{G(T_{[0,t(x))}(\Phi),C)}(x)$
is finite by the fact that $T_{[0,t(x))}(\Phi)$ is an homogeneous
Poisson p.p. of intensity smaller than
$n\lambda^{0}$ and the induction hypothesis. If $t(x)\geq \frac{n\lambda_{0}}{\lambda}$,
we rewrite $K^{*}_{G(\Phi,C)}(x)$ as:
\begin{align*}
% &K^{*}_{G(\Phi,C)}(x)=\\
 & K^{*}_{G\left(T_{[\frac{n\lambda_{0}}{\lambda},1]}(\Phi),C\right)}(x)\bigcup
\left(
\bigcup_{y \in K^{*}_{G\left(T_{[\frac{n\lambda_{0}}{\lambda},1]}(\Phi),C\right)}(x)}
K^{*}_{G\left(T_{[0,\frac{n\lambda_{0}}{\lambda})}(\Phi)\cup \left(y\right)^{+},C\right)}(y)\right)
.
\end{align*}
The first term is the set of points accessible from $x$ with timer larger than
$\frac{n\lambda_{0}}{\lambda}$; the second term is the set of points 
accessible from $x$ with timer smaller than $\frac{n\lambda_{0}}{\lambda}$. The second term
can in turn be decomposed as the union of the sets of points accessible from $y$ with timer
smaller than $\frac{n\lambda_{0}}{\lambda}$,  for all $y$
accessible form $x$ with timer larger than $\frac{n\lambda_{0}}{\lambda}$.
Since $G\left(T_{[\frac{n\lambda_{0}}{\lambda},1]}(\Phi),C\right)$ has no
infinite positive percolation, we have that the first term is finite and the 
union in the second term is a finite union. Moreover, by Lemma \ref{L:Palm finite per cond}
and the induction hypothesis, we have that each term in this finite union
is finite. Thus $K^{*}_{G(\Phi,C)}(x)$ is finite, which terminates the proof.
 \begin{flushright}
 $\Box$
\end{flushright} 

  \subsubsection{Sufficient Condition}
  Let us first state the main result:
  \begin{proposition}\label{prop1}
  If the contention field $C$ satisfies (\ref{E:CSMA well def cond}), then the $\infty$-Mat\'ern model of the Poisson point process is well-defined in the sense that
its thinning indicator $e_{\infty}(.)$ is well defined.
\end{proposition}
The proof of this proposition is based on the following lemma:

\begin{lemma}\label{L:finite per CSMA well-def}
For all p.p.s $\Phi$ such that $G(\Phi,C)$ has no infinite
positive percolation, we have $e_{\infty}(x)= e(K^{*}_{G(\Phi,C)}(x),x)$, where $e(.)$
is the thinning indicator function defined in (\ref{E:Thin ind func}). This is also the only
solution of Equation (\ref{E:CSMA thin ind}).
\end{lemma}

\emph{Proof.}
Since $G(\Phi,C)$ has no infinite positive percolation, $e(K^{*}_{G(\Phi,C)}(x),x)$ is well defined for all $x$. As 
$K^{*}_{G(\Phi,C)}(x)=x\cup\left(\bigcup_{y \textnormal{ s.t. }(x,y)\in G(\Phi,C)}K^{*}_{G(\Phi,C)}(y)\right),$ we have that $e(K^{*}_{G(\Phi,C)}(x),x)$
satisfies Equation (\ref{E:CSMA thin ind}) by definition.\\
For uniqueness, let $e(x)$ be a solution of (\ref{E:CSMA thin ind}).
We prove by induction on the size of $K^{*}_{G(\Phi,C)}(x)$ that 
$e(x)= e(K^{*}_{G(\Phi,C)}(x),x)$.\\
For any $x$ such that $|K^{*}_{G(\Phi,C)}(x)|=1$, we have $K^{*}_{G(\Phi,C)}(x)=\{x\}$ or equivalently, $x$ has out degree $0$ in the conflict graph.
Hence \[e(x)=e(K^{*}_{G(\Phi,C)}(x),x)=1.\]
 Note that such $x$ always exists by the 
assumption that $G(\Phi,C)$ has no infinite percolation.\\
Suppose that for all $y$ such that $|K^{*}_{G(\Phi,C)}(y)|< n$,
we have $e(y)=e(K^{*}_{G(\Phi,C)}$ $(y),y)$. Let us consider an $x$
such that $$|K^{*}_{G(\Phi,C)}(x)|= n$$ (if there are any).
For $y$ such that $(x,y)\in G(\Phi,C)$, we have $|K^{*}_{G(\Phi,C)}(y)|< n$.
Thus, $e(y)= e(K^{*}_{G(\Phi,C)}(y),y)$.
Then:
\begin{align*}
 e(x)&= \prod_{y \in T_{[0,t(x))}(\Phi)}(1-e(y)C(x,y))= \prod_{y \in \Phi \textnormal{ s.t. } t(y)<t(x)}(1-e(y)C(x,y))\\
     &= \prod_{y \in \Phi \textnormal{ s.t. } t(y)<t(x), C(x,y)=1}(1-e(y))= \prod_{y \in \Phi \textnormal{ s.t. } (x,y)\in G(\Phi,C)}(1-e(y))\\
     &= \prod_{y \in \Phi \textnormal{ s.t. } (x,y)\in G(\Phi,C)}(1-e(K^{*}_{G(\Phi,C)}(y),y))= e(K^{*}_{G(\Phi,C)}(x),x).
\end{align*}
Thus, for all $x$ such that $|K^{*}_{G(\Phi,C)}(x)|<n+1$ we have
 \[e(x)=e(K^{*}_{G(\Phi,C)}(x),x).\]
This completes the proof of uniqueness.
 \begin{flushright}
 $\Box$
\end{flushright}

Now, it is straightforward to see that when
(\ref{E:CSMA well def cond}) is satisfied, the conflict graph associated with $\Phi$
has no infinite positive percolation (Lemma \ref{L:large finite per cond}). Thus, the $\infty$-Mat\'ern thinning
indicators are well-defined (Lemma \ref{L:finite per CSMA well-def}) and are the only solution of (\ref{E:CSMA thin ind}). 
Note that the arguments in the above proof can be extended to the case where $\Phi$ is inhomogeneous to obtain the following proposition.
 \begin{proposition}
  Let $\Phi$ be an inhomogeneous Poisson p.p. of intensity measure $\Lambda$. If we have:
  \begin{align}
  \sup_{x \in \mathbb{R}^{2}}\int_{\mathbb{R}^{2}}h(x,y)d\Lambda(y)= \mathcal{N} <\infty, \label{E:inhom CSMA well def cond}
  \end{align}
  then the $\infty$-Mat\'ern model is well-defined in the sense that its thinning indicators $e_{\infty}(x)$ are well defined.
 \end{proposition}
 \subsubsection{The $\infty$-Mat\'ern Model as Limit of the $k$-Mat\'ern Models}
In this part we explain the meaning of the subscript $\infty$ in the definition of the $\infty$-Mat\'ern model. More precisely, we prove that if the conflict graph has no infinite positive percolation, then the thinnning indicators are such that:
\begin{align*}
 \forall x\in \Phi \qquad e_{\infty}(x)= \lim_{k \rightarrow \infty}e_{k}(x) \textnormal{ a.s. }
\end{align*}
Recall that $e(.)$ is the thinning indicator function defined in (\ref{E:Thin ind func}) and $K^{k}(x)$ is the $k$-step positive connected component of $x$. Let us begin with an important result relating the thinning indicators of the $k$-Mat\'ern model with the conflict indicators.
\begin{proposition}\label{P: k Matern char}
 For all $k \in \mathbb{N}$: $e_{k}(x)= e(K^{k}_{G(\Phi,C)}(x),x)$.
\end{proposition}
\emph{Proof}\\
  We prove this by induction. The base case $k=0$ is trivial since $e_{0}(x)=e(K^{0}_{G(\Phi,C)}(x),x)=1.$
For $k\geq 1$ we have by definition:
\begin{align*}
 e_{k}(x)&= \prod_{y \in T_{[0,t(x))}(\Phi)}(1-e_{k-1}(y)C(x,y))\\
	 &= \prod_{y \in \Phi \textnormal{ s.t. } t(y)<t(x),C(x,y)=1}(1-e_{k-1}(y))\\
	 &= \prod_{y \in \Phi \textnormal{ s.t. } (x,y) \in G(\Phi,C)}(1-e_{k-1}(y)).
\end{align*}
Now, $e_{k-1}(y)=e(K^{k-1}_{G(\Phi,C)}(y),y)$ by induction hypothesis. We have:
\begin{align*}
 e_{k}(x)= \prod_{y \in \Phi \textnormal{ s.t. } (x,y) \in G(\Phi,C)}(1-e(K^{k-1}_{G(\Phi,C)}(y),y))= e(K^{k}_{G(\Phi,C)}(x),x).
\end{align*}
\begin{flushright}
 $\Box$
\end{flushright}
\begin{proposition}\label{P: k Matern convergence to CSMA}
 If the conflict graph $G(\Phi,C)$ has no infinite positive percolation, then:
\[ \forall x \in \Phi \qquad e_{\infty}(x)= \lim_{k \rightarrow \infty}e_{k}(x) \textnormal{ a.s. }\]
\end{proposition}
\emph{Proof}\\
If the conflict graph has no infinite positive percolation, then $K^{*}_{G(\Phi,C)}(x)$ is finite and  $K^{*}_{G(\Phi,C)}(x)=K^{|K^{*}_{G(\Phi,C)}(x)|}_{G(\Phi,C)}(x)$ for all $x$.
It then follows from Lemma \ref{L:finite per CSMA well-def} and Proposition \ref{P: k Matern char} that the sequence $e_{k}(x)$ is constant for $k>|K^{*}_{G(\Phi,C)}(x)|$ and this constant is equal to 
$e_{\infty}(x)$.
\begin{flushright}
 $\Box$
\end{flushright}
\section{The Generating Functional}\label{S:Anal}
\subsection{Some Background}\label{SS:Anal back}
For a p.p. $\Xi$, for all functions $v$ from $\mathbb{R}^{2}$ to $\mathbb{R}^{+}$ satisfying:
\begin{align}
 \textbf{E}\left[\sum_{x \in \Xi}|1-v(x)|\right]<\infty, \label{E: Gen def con}
\end{align}
we have $\sum_{x \in \Xi}|1-v(x)|<\infty$ a.s. Since $|\log(t)|<|1-t|$ for all $t>0$, we have that $\sum_{x \in \Xi}|\log(v(x))|<\infty$ a.s.. Hence the sum  $\sum_{x \in \Xi}\log(v(x))$ is well-defined a.s. and is itself a random variable. Thus, we can define:
\begin{align*}
 G_{\Xi}(v)= \textbf{E}\left[\prod_{x \in \Xi}v(x)\right]=\textbf{E}\left[\exp\left(\sum_{x \in \Xi}\log(v(x))\right)\right],
\end{align*}
as the generating functional of $\Xi$. It encodes all the information about the distribution of $\Xi$. Many properties of $\Xi$ such as the void probability, the contact distribution, the laplace transform of its Shot-Noise process can be computed in terms of its generating functional (see \cite{Stoyan,DarVer1}).\\ 
In this section we are interested in the generating functionals
$G_{M^{k}(\Phi,C)}(v(.))$, where $\Phi$ is a Poisson p.p. of intensity measure $\Lambda$
and $v$ is a function from $\mathbb{R}^{2}$ to $[0,1]$ such that:
\begin{align}
 \int_{\mathbb{R}^{2}}|1-v(x)|\Lambda(dx) <\infty. \label{E: Pois gen func cond}
\end{align}
\emph{Remark}:
\begin{itemize}
 \item The condition (\ref{E: Pois gen func cond}) is equivalent to
$\textbf{E}[\sum_{x \in \Phi}|1-v(x)|]<\infty$ (Campbell's formula). Thus, this condition guarantees
that the generating functional of $\Phi$ is 
well-defined at $v$. Since $M^{k}(\Phi,C)$ is dominated by $\Phi$,
the generating functional of $M^{k}(\Phi,C)$ is also well-defined at $v$.
\item Although there are functions $v$ such that the generating functional of $M^{k}(\Phi,C)$
is well-defined at $v$ while  the generating functional of $\Phi$ is not, we will not consider 
this case here for two reasons:
\begin{itemize}
\item All informations about the distribution of a p.p. are encoded in the generating
functional with $v$ having bounded support. These functions $v$ satisfy
(\ref{E: Pois gen func cond}) automatically.
       %\item For applications requires computing the generating functional of the p.p. of interested with some specific $v$. The function $v$ is always well-behaved, i.e. satisfies (\ref{E: Pois gen func cond}).
\item As we will see later, this condition guarantees nice convergence properties.
\end{itemize}
\end{itemize}
In order to study the functionals $G_{M^{k}(\Phi,C)}(v(.))$,
it will be useful to extend our model to the case where the timers are no longer restricted to $[0,1]$
and take values ``uniformly'' on $\mathbb{R}^{+}$. This makes our p.p. $\Phi$ no 
longer a marked Poisson p.p. but a Poisson p.p. on $\mathbb{R}^{2}\times\mathbb{R}^{+}$,
with the last dimension being the timer. The intensity measure of 
$\Phi$ satisfies
$\textbf{E}[\Phi(A\times[0,t])]=t\Lambda(A),$
for all Borel sets $A$ and all positive real numbers $t$.
The description of the timer based thinning and the $\Psi$ thinning are kept unchanged.
Our aim is to study the following functionals:
\begin{align}
f_{k,\Lambda}(t,v):= G_{M^{k}(T_{[0,t)}(\Phi),C)}(v(.)). \label{E: F def}
\end{align}
The link between the model with timers restricted to $[0,1]$ and the last model
is that the generating functional of the former is $f_{k,\Lambda}(1,v)$.

%In the above definition, if $\Lambda$ is the Lebesgue measure then the subscript is omitted.
In this section, we gather our major results regarding the evolution over time $t$ of the $f_{k,\Lambda}(t,v)$ functional.
More precisely, we prove that it is the solution of some differential equation w.r.t. $t$.
First we show that it is continuous.  
\begin{proposition}\label{P:f cont}
For all $k$ and all $v$ satisfying the condition (\ref{E: Pois gen func cond}),
$f_{k,\Lambda}(t,v)$ is continuous in $t$.
\end{proposition}
\emph{Proof.}
For all $t_{1}<t_{2}$ and all $k \in \mathbb{N}$:
\[T_{[0,t_{1})}(M^{k}(T_{[0,t_{2})}(\Phi),C))=M^{k}(T_{[0,t_{1})}(\Phi),C). \]
Thus, for all $t$ positive and $\epsilon$ positive:
\begin{align*}
& f_{k,\Lambda}(t+\epsilon,v)\\
&=\textbf{E}\left[\prod_{x \in M^{k}(T_{[0,t+\epsilon)}(\Phi),C)}v(x)\right]\\ 
&=\textbf{E}\left[\prod_{x \in T_{[0,t)}(M^{k}(T_{[0,t+\epsilon)}(\Phi),C))}v(x)\prod_{x \in T_{[t,t+\epsilon)}(M^{k}(T_{[0,t+\epsilon)}(\Phi),C))}v(x)\right]\\
&=\textbf{E}\left[\prod_{x \in M^{k}(T_{[0,t)}(\Phi),C)}v(x)\prod_{x \in T_{[t,t+\epsilon)}(M^{k}(T_{[0,t+\epsilon)}(\Phi),C))}v(x)\right].
\end{align*}
Since $M^{k}(T_{[0,t+\epsilon)}(\Phi),C)$ is dominated by
$T_{[0,t+\epsilon)}(\Phi)$, $T_{[t,t+\epsilon)}(M^{k}(T_{[0,t+\epsilon)}(\Phi),C))$
is dominated by $T_{[t,t+\epsilon)}(T_{[0,t+\epsilon)}(\Phi))=T_{[t,t+\epsilon)}(\Phi)$. Let $v_{+}(x)=\max\{1,v(x)\}$ and $v_{-}(x)=\min\{1,v(x)\}$, we have a.s.
 \begin{align*}
&\prod_{x \in T_{[t,t+\epsilon)}(\Phi)}v_{+}(x)\geq \prod_{x \in T_{[t,t+\epsilon)}(M^{k}(T_{[0,t+\epsilon)}(\Phi),C))}v_{+}(x) \\
&\geq \prod_{x \in T_{[t,t+\epsilon)}(M^{k}(T_{[0,t+\epsilon)}(\Phi),C))}v(x)\\
&\geq \prod_{x \in T_{[t,t+\epsilon)}(M^{k}(T_{[0,t+\epsilon)}(\Phi),C))}v_{-}(x) \geq \prod_{x \in T_{[t,t+\epsilon)}(\Phi)}v_{-}(x).
\end{align*}
Thus:
 \begin{align*}
&f_{k,\Lambda}(t,v)\exp\{\epsilon\int_{\mathbb{R}^{2}}|1-v(x)|\Lambda(dx)\}\geq f_{k,\Lambda}(t,v)\exp\{\epsilon\int_{\mathbb{R}^{2}}(v_{+}(x)-1\Lambda(dx)\}\\ 
&= \textbf{E}\left[\prod_{x \in M^{k}(T_{[0,t)}(\Phi),C)}v(x)\prod_{x \in T_{[t,t+\epsilon)}(\Phi)}v_{+}(x)\right]\geq \textbf{E}\left[\prod_{x \in M^{k}(T_{[0,t+\epsilon)}(\Phi),C)}v(x)\right]\\
&= f_{k,\Lambda}(t+\epsilon,v)\geq \textbf{E}\left[\prod_{x \in M^{k}(T_{[0,t)}(\Phi),C)}v(x)\prod_{x \in T_{[t,t+\epsilon)}(\Phi)}v(x)\right]\\
&= \textbf{E}\left[\prod_{x \in M^{k}(T_{[0,t)}(\Phi),C)}v(x)\prod_{x \in T_{[t,t+\epsilon)}(\Phi)}v_{-}(x)\right]\\
&= f_{k,\Lambda}(t,v)\exp\left\{-\epsilon\int_{\mathbb{R}^{2}}(1-v_{-}(x))\Lambda(dx)\right\}\\
&\geq f_{k,\Lambda}(t,v)\exp\left\{-\epsilon\int_{\mathbb{R}^{2}}|1-v(x)|\Lambda(dx)\right\}.
\end{align*}
Doing the same, we have:
\begin{align*}
&f_{k,\Lambda}(t-\epsilon,v)\exp\{\epsilon\int_{\mathbb{R}^{2}}|1-v(x)|\Lambda(dx)\}\geq f_{k,\Lambda}(t,v)\\
&\geq f_{k,\Lambda}(t-\epsilon,v)\exp\left\{-\epsilon\int_{\mathbb{R}^{2}}|1-v(x)|\Lambda(dx)\right\}.
\end{align*}

Letting $\epsilon$ go to $0$ completes our proof.
\begin{flushright}
 $\Box$.
\end{flushright}
\subsection{Differential Equation for the $\infty$-Mat\'ern Model}
Let $H$ be the map that associates to a function $v$ from ${\mathbb{R}^{2}}\to [0,1]$ and a point
$x\in \mathbb{R}^{2}$ the function 
\begin{align}
H(v,x)(y)= v(y)(1-h(x,y)). \label{E:h transform}
\end{align}
It is easy to see that $H(v,x)(.)$ is a function from ${\mathbb{R}^{2}}$ to $[0,1]$.
\begin{theorem}\label{T:CSMA eq}
For any locally finite measure $\Lambda$, the functional $f_{\infty,\Lambda}$ satisfies the following equation:
\begin{eqnarray}
\label{E:CSMA eq}
 f_{\infty,\Lambda}(0,v(.))&=&1\\
 \frac{df_{\infty,\Lambda}(t,v(.))}{dt}&=& -\int_{\mathbb{R}^{2}}f_{\infty,\Lambda}(t,H(v,x))(1-v(x))\Lambda(dx).
\nonumber
\end{eqnarray}
\end{theorem}
The main idea behind this theorem is to divide the p.p.
$\Phi$ into thin layers $T_{[t,t+\epsilon)}(\Phi)$ such that the points in
each layer are so sparse that the effect of contention is negligible. Then, one can consider the retained p.p.s in these layers as Poisson p.p.s \\
In particular, we need to prove that:
\begin{align}
%& &
\lim_{\epsilon \rightarrow 0}\frac{f(t+\epsilon,v)-f(t,v)}{\epsilon} =
-\int_{\mathbb{R}^{2}}f(t,H(v,x))(1-v(x))\Lambda(dx) \label{E:left limit}\\
%& &
\lim_{\epsilon \rightarrow 0} \frac{f(t,v)-f(t-\epsilon,v)}{\epsilon}= 
-\int\limits_{\mathbb{R}^{2}}f(t,H(v,x))(1-v(x))\Lambda(dx)\label{E:right limit},
\end{align}
for any locally finite $\Lambda$.
 We will need the following lemmas:

\begin{lemma}\label{L: main eq}
 For any p.p. $\Xi$ and any function $v$ taking value in $\mathbb{R}^{+}$ such that:
\begin{align*}
 \textbf{E}\left[\sum_{x \in \Xi}|1-v(X)|\right]<\infty,
\end{align*}
we have:
\begin{align}
 \prod_{x  \in \Xi}v(x)= 1+\sum_{i=1}^{\infty}(-1)^{i}\sum_{(x_{1},\cdots,x_{i}) \in \Xi^{i!}}\prod_{j=1}^{i}(1-v(x_{j})) \textnormal{ a.s.},\label{E: main eq}
\end{align}
where $\Xi^{i!}$ is the set of unordered $i$-tuples of points in $\Xi$.
\end{lemma}
\emph{Proof.}
By the same argument as those at the beginning of the section, we have $\prod_{x  \in \Xi}v(x)= \exp\{\sum_{x \in \Xi}\log(v(x))\}$ is well-defined and a.s. finite.\\
Now, to prove that the series in the right hand side of (\ref{E: main eq}) converges, it is enough to show that:
\begin{align*}
  1+\sum_{i=1}^{\infty}\sum_{(x_{1},\cdots,x_{i}) \in \Xi^{i!}}\prod_{j=1}^{i}|1-v(x_{j})|<\infty \textnormal{ a.s.}.
\end{align*}
 Note that:
\begin{align*}
& \sum_{(x_{1},\cdots,x_{i}) \in \Xi^{i!}}\prod_{j=1}^{i}|1-v(x_{j})|=\frac{1}{i!}\sum_{x_{1},\cdots,x_{i} \textnormal{ mutually different} \in \Xi}\textnormal{ }\prod_{j=1}^{i}|1-v(x_{j})|\\
&\leq \frac{1}{i!}\left(\sum_{x \in \Xi}|1-v(x)|\right)^{i}.
\end{align*}
Hence:
\begin{align*}
&  1+\sum_{i=1}^{\infty}\sum_{(x_{1},\cdots,x_{i}) \in \Xi^{i!}}\prod_{j=1}^{i}|1-v(x_{j})|\leq 1+\sum_{i=1}^{\infty}\frac{1}{i!}\left(\sum_{x \in \Xi}|1-v(x)|\right)^{i}\\
&=e^{\left(\sum_{x \in \Xi}|1-v(x)|\right)}<\infty.
\end{align*}
Now the equality can be obtained by writing:
\begin{align*}
 \prod_{x\in \Xi}v(x)= \prod_{x \in \Xi}(1-(1-v(x))).
\end{align*}
\begin{flushright}
 $\Box$
\end{flushright}

\begin{lemma}\label{L:limit 0}
 For any $k$, including $\infty$ and any locally finite measure $\Lambda$, we have:
\begin{align*}
 &\int_{\mathbb{R}^{2}}\textbf{1}_{v(x)<1}(1-v(x))\Lambda(dx)+\int_{\mathbb{R}^{2}}\textbf{1}_{v(x)\geq 1}(1-v(x))e^{-\epsilon\int_{\mathbb{R}^{2}}(1-h(x,y))\Lambda(dy)}\Lambda(dx)\\
&+\frac{e^{\epsilon\int_{\mathbb{R}^{2}}|1-v(x)|\Lambda(dx)}-1-\epsilon\int_{\mathbb{R}^{2}}|1-v(x)|\Lambda(dx)}{\epsilon}\geq \frac{1-f_{k,\Lambda}(\epsilon,v)}{\epsilon}\\
 &\geq \int_{\mathbb{R}^{2}}\textbf{1}_{v(x)>1}(1-v(x))\Lambda(dx)+\int_{\mathbb{R}^{2}}\textbf{1}_{v(x)\leq 1}(1-v(x))e^{-\epsilon\int_{\mathbb{R}^{2}}(1-h(x,y))\Lambda(dy)}\Lambda(dx)\\
&-\frac{e^{\epsilon\int_{\mathbb{R}^{2}}|1-v(x)|\Lambda(dx)}-1-\epsilon\int_{\mathbb{R}^{2}}|1-v(x)|\Lambda(dx)}{\epsilon}.
\end{align*}
In particular:
\begin{align}
 \lim_{\epsilon \rightarrow 0}\frac{f_{k,\Lambda}(\epsilon,v)-1}{\epsilon}=-\int_{\mathbb{R}^{2}}(1-v(x))\Lambda(dx). \label{E:quarsi-Poisson}
\end{align}
\end{lemma}
\emph{Proof}.
We first use Lemma \ref{L: main eq} with $M^{k}(T_{[0,\epsilon)}(\Phi),C)$ in place of $\Xi$ to get:
\begin{align*}
 1-f_{k,\Lambda}(\epsilon,v)=\sum_{i=1}^{\infty}(-1)^{i+1}\textbf{E}\left[\sum_{(x_{1},\cdots,x_{i}) \in M^{k}(T_{[0,\epsilon)}(\Phi),C)^{i!}}\prod_{j=1}^{i}(1-v(x_{j}))\right]
\end{align*}
For any $i>1$, we have:
\begin{align*}
& \left|\textbf{E}\left[\sum_{(x_{1},\cdots,x_{i}) \in M^{k}(T_{[0,\epsilon)}(\Phi),C)^{i!}}\prod_{j=1}^{i}(1-v(x_{j}))\right]\right|\\
&= \frac{1}{i!}\left|\textbf{E}\left[\sum_{x_{1},\cdots,x_{i} \textnormal{ mutually different} \in  M^{k}(T_{[0,\epsilon)}(\Phi),C)^{i!}}\prod_{j=1}^{i}(1-v(x_{j}))\right]\right|\\
&= \frac{1}{i!}\int_{(\mathbb{R}^{2})^{i}}\prod_{j=1}^{i}\left|1-v(x_{j})\right|m^{i}_{k,\Lambda}(dx_{1},\cdots,dx_{j}),
\end{align*}
where $m^{i}_{k,\Lambda}$ is the $i$-factorial moment measure of $M^{k}(T_{[0,\epsilon)}(\Phi),C)$. Now, by the fact that $M^{k}(T_{[0,\epsilon)}(\Phi),C)$ is dominated by $T_{[0,\epsilon)}(\Phi)$, its $i$-factorial moment measure is dominated by that of $T_{[0,\epsilon)}(\Phi)$, which is known to be $\epsilon^{i}\prod_{j=1}^{i}\Lambda(dx_{j})$. Thus,
\begin{align*}
 & \left|\textbf{E}\left[\sum_{(x_{1},\cdots,x_{i}) \in M^{k}(T_{[0,\epsilon)}(\Phi),C)^{i!}}\prod_{j=1}^{i}(1-v(x_{j}))\right]\right|\\
&\leq \frac{1}{i!}\epsilon^{i}\int_{(\mathbb{R}^{2})^{i}}\prod_{j=1}^{i}\left|1-v(x_{j})\right|\prod_{j=1}^{i}\Lambda(dx_{j}) =\frac{1}{i!}\epsilon^{i}\left(\int_{\mathbb{R}^{2}}|1-v(x)|\Lambda(dx)\right)^{i}.
\end{align*}
Hence:
\begin{align*}
&-\frac{e^{\epsilon\int_{\mathbb{R}^{2}}|1-v(x)|\lambda(dx)}-1-\epsilon\int_{\mathbb{R}^{2}}|1-v(x)|\lambda(dx)}{\epsilon}= -\sum_{i=2}^{\infty} \frac{1}{i!}\epsilon^{i}\left(\int_{\mathbb{R}^{2}}|1-v(x)|\Lambda(dx)\right)^{i}\\
& \leq \sum_{i=2}^{\infty}(-1)^{i+1}\textbf{E}\left[\sum_{(x_{1},\cdots,x_{i}) \in M^{k}(T_{[0,\epsilon)}(\Phi),C)^{i!}}\prod_{j=1}^{i}(1-v(x_{j}))\right]\\
& \leq \sum_{i=2}^{\infty} \frac{1}{i!}\epsilon^{i}\left(\int_{\mathbb{R}^{2}}|1-v(x)|\Lambda(dx)\right)^{i}=\frac{e^{\epsilon\int_{\mathbb{R}^{2}}|1-v(x)|\lambda(dx)}-1-\epsilon\int_{\mathbb{R}^{2}}|1-v(x)|\lambda(dx)}{\epsilon}.
\end{align*}
Now we bound $\textbf{E}\left[\sum_{x \in M^{k}(T_{[0,\epsilon)}(\Phi),C)}(1-v(x))\right]$. Let us define the p.p. 
\[\varPi_{\epsilon}=\{x \in T_{[0,\epsilon)}(\Phi) \textnormal{ s.t. } C(x,y)=0 \textnormal{ } \forall y \in T_{[0,\epsilon)}(\Phi)\}.\]
This is the process of points of $T_{[0,\epsilon}(\Phi)$ that do not have any contender. Its first moment measure can be computed as:
\begin{align*}
\textbf{E}[\varPi_{\epsilon}(B)]&=\textbf{E}\left[\sum_{x \in T_{[0,\epsilon}(\Phi)}\prod_{y \in T_{[0,\epsilon}(\Phi), y \neq x}(1-C(x,y))\right]\\
&=\int_{B}e^{-\epsilon\int_{\mathbb{R}^{2}}(1-h(x,y))\Lambda(dy)}\Lambda(dx).
\end{align*}
We can easily see that $ M^{k}(T_{[0,\epsilon)}(\Phi),C)$ is dominated by $T_{[0,\epsilon)}(\Phi)$ and dominates $\varPi_{\epsilon}$. Hence:
\begin{align*}
 &\textbf{E}\left[\sum_{x \in M^{k}(T_{[0,\epsilon)}(\Phi),C)}(1-v(x))\right]\\
 &=\textbf{E}\left[\sum_{x \in M^{k}(T_{[0,\epsilon)}(\Phi),C)}\textbf{1}_{v(x)>1}(1-v(x))\right]+\textbf{E}\left[\sum_{x \in M^{k}(T_{[0,\epsilon)}(\Phi),C)}\textbf{1}_{v(x)\leq 1}(1-v(x))\right]\\
 &\geq \textbf{E}\left[\sum_{x \in T_{[0,\epsilon)}(\Phi)}\textbf{1}_{v(x)>1}(1-v(x))\right]+\textbf{E}\left[\sum_{x \in \varPi_{\epsilon}}\textbf{1}_{v(x)\leq 1}(1-v(x))\right]\\
 &=\int_{\mathbb{R}^{2}}\textbf{1}_{v(x)>1}(1-v(x))\Lambda(dx)+\int_{\mathbb{R}^{2}}\textbf{1}_{v(x)\leq 1}(1-v(x))e^{-\epsilon\int_{\mathbb{R}^{2}}(1-h(x,y))\Lambda(dy)}\Lambda(dx).
\end{align*}
Doing similarly, we have:
\begin{align*}
 &\textbf{E}\left[\sum_{x \in M^{k}(T_{[0,\epsilon)}(\Phi),C)}(1-v(x))\right]\\
 &\leq \int_{\mathbb{R}^{2}}\textbf{1}_{v(x)\leq 1}(1-v(x))\Lambda(dx)+\int_{\mathbb{R}^{2}}\textbf{1}_{v(x)> 1}(1-v(x))e^{-\epsilon\int_{\mathbb{R}^{2}}(1-h(x,y))\Lambda(dy)}\Lambda(dx).
\end{align*}
Combining the bound for $\textbf{E}\left[\sum_{x \in M^{k}(T_{[0,\epsilon)}(\Phi),C)}(1-v(x))\right]$ and the bound for $\sum_{i=2}^{\infty}(-1)^{i+1}\textbf{E}\left[\sum_{(x_{1},\cdots,x_{i}) \in M^{k}(T_{[0,\epsilon)}(\Phi),C)^{i!}}\prod_{j=1}^{i}(1-v(x_{j}))\right]$, we get the desired result.
\begin{flushright}
 $\Box$
\end{flushright}
At this point, it is useful to compare the result just obtained in the above lemma with that of a Poisson p.p. Note that the generating functional of $T_{[0,\epsilon)}(\Phi)$, which is a Poisson p.p. of intensity measure $\epsilon\Lambda$, is:
\begin{align*}
 G_{T_{[0,\epsilon)}(\Phi)}(v)=e^{-\epsilon\int_{\mathbb{R}^{2}}(1-v(x))\Lambda(dx)}.
\end{align*}
Hence:
\begin{align*}
 \lim_{\epsilon \rightarrow 0}\frac{G_{T_{[0,\epsilon)}(\Phi)}(v)-1}{\epsilon}=-\int_{\mathbb{R}^{2}}(1-v(x))\Lambda(dx).
\end{align*}
Thus, Lemma \ref{L:limit 0} justifies our intuition that when the time scale is small, the effect of contention is negligible, and we can consider the thin layer of retained points as a Poisson p.p. Such a property, which we call the \emph{quasi-Poisson} property, plays an important role in the subsequent study in this paper.

Now we can proceed to the proof of Theorem \ref{T:CSMA eq}. 
Note that for $t>s$:
\begin{eqnarray*}
 f(t,v)-f(s,v) = \textbf{E}\left[\prod_{x \in M^{\infty}(T_{[0,s)}(\Phi))}v(x)
\left(\prod_{y \in T_{[s,t)}(M^{\infty}(T_{[0,t)}(\Phi)))}v(y)-1\right)\right].
\end{eqnarray*}
In order to evaluate the last expression, we need the following conditional probability:
\begin{align*}
 \textbf{E}\left[\left.\prod_{y \in T_{[s,t)}(M^{\infty}(T_{[0,t)}(\Phi)))}v(y)\right|T_{[0,s)}(\Phi)\right].
\end{align*}
Let $A$ be any countable set of points such that $\prod_{x \in A}(1-h(z,x))$ is well defined for all $x$ in $\mathbb{R}^{2}$. For all measures $\Lambda$ and all denumerable sets $A$ we define the $A$-constrained measure $\Lambda_{A}$ as the unique measure such that $\Lambda_{A}(dz)=\prod_{x \in A}(1-h(z,x))\Lambda(dz)$. This is the intensity of the process of points belonging to a Poisson p.p. of intensity measure $\Lambda$ which do not contend with any point in $A$.
 Note that $M^{\infty}(T_{[0,s)}(\Phi),C)$ is always such that $\prod_{z \in M^{\infty}(T_{[0,s)}(\Phi),C)}(1-h(z,x))$ is well defined for all $x$ in $\mathbb{R}^{2}$. We have:  
\begin{align*}
 \textbf{E}\left[\left.\prod_{y \in T_{[s,t)}(M^{\infty}(T_{[0,t)}(\Phi)))}v(y)\right|T_{[0,s)}(\Phi)\right]= f_{\infty,\Lambda_{M^{\infty}(T_{[0,t)}(\Phi)}}(s-t,v).
\end{align*} 
Put $s=t$ and $t=t+\epsilon$. By using the limit in Lemma \ref{L:limit 0}, we have:  
\begin{align*}
%& &
&\lim_{\epsilon\rightarrow 0}\frac{f(t+\epsilon,v)-f(t,v)}{\epsilon}\\
&=\lim_{\epsilon\rightarrow 0}\textbf{E}\left[\prod_{z \in M^{\infty}(T_{[0,t)}(\Phi),C)}v(z)
\frac{f_{\infty,\Lambda_{M^{\infty}(T_{[0,t)}(\Phi),C)}}(\epsilon,v)-1}{\epsilon}\right]\\
&=\textbf{E}\left[\prod_{z \in M^{\infty}(T_{[0,t)}(\Phi),C)}v(z)
\lim_{\epsilon\rightarrow 0}\frac{f_{\infty,\Lambda_{M^{\infty}(T_{[0,t)}(\Phi),C)}}(\epsilon,v)-1}{\epsilon}\right]\\
&=-\textbf{E}\left[\prod_{z \in M^{\infty}(T_{[0,t)}(\Phi),C)}v(z)
\int_{\mathbb{R}^{2}}(1-v(x))\prod_{z \in M^{\infty}(T_{[0,s)}(\Phi),C)}(1-h(x,z))\Lambda(dx)\right]\\
&=-\int_{\mathbb{R}^{2}}\textbf{E}\left[\prod_{z \in M^{\infty}(T_{[0,t)}(\Phi),C)}v(z)(1-h(x,z))\right](1-v(x))\Lambda(dx)\\
&=-\int_{\mathbb{R}^{2}}f_{\infty,\Lambda}(t,H(v,x))(1-v(x))\Lambda(dx).
\end{align*}
Now we put $s=t-\epsilon$ and $t=t$. Using the lower bound in Lemma \ref{L:limit 0}, we have:
\begin{align*}
&\frac{f(t,v)-f(t-\epsilon,v)}{\epsilon}\\
&=\textbf{E}\left[\prod_{z \in M^{\infty}(T_{[0,t-\epsilon)}(\Phi),C)}v(z)
\frac{f_{\infty,\Lambda_{M^{\infty}(T_{[0,t-\epsilon)}(\Phi),C)}}(\epsilon,v)-1}{\epsilon}\right]\\
&=\textbf{E}\left[\prod_{z \in M^{\infty}(T_{[0,t)}(\Phi),C)}v(z)
\frac{f_{\infty,\Lambda_{M^{\infty}(T_{[0,t-\epsilon)}(\Phi),C)}}(\epsilon,v)-1}{\epsilon}\right]\\
&\geq-\textbf{E}\Biggl[\prod_{z \in M^{\infty}(T_{[0,t-\epsilon)}(\Phi),C)}v(z)
\Biggl(\int_{\mathbb{R}^{2}}\textbf{1}_{v(x)\leq 1}(1-v(x))\prod_{z \in M^{\infty}(T_{[0,t-\epsilon)}(\Phi),C)}\\
&(1-h(x,z))\Lambda(dx)+\int_{\mathbb{R}^{2}}\textbf{1}_{v(x)>1}(1-v(x))\prod_{z \in M^{\infty}(T_{[0,t-\epsilon)}(\Phi),C)}(1-h(x,z))\\
&e^{-\epsilon\int_{\mathbb{R}^{2}}(1-h(x,y))\prod_{z \in M^{\infty}(T_{[0,t-\epsilon)}(\Phi),C)}(1-h(y,z))\Lambda(dy)}\Lambda(dx)+\\
& \vartheta\Biggl(\epsilon,\int_{\mathbb{R}^{2}}(1-v(x))\prod_{z \in M^{\infty}(T_{[0,t-\epsilon)}(\Phi),C)}(1-h(x,z))\Lambda(dx)\Biggr)\Biggr)\Biggr].
\end{align*}
where $\vartheta(x,y)=\frac{e^{xy}-1-xy}{x}$.
Then, by the upper bound:
\begin{align*}
&\frac{f(t,v)-f(t-\epsilon,v)}{\epsilon}\\
&=\textbf{E}\left[\prod_{z \in M^{\infty}(T_{[0,t-\epsilon)}(\Phi),C)}v(z)
\frac{f_{\infty,\Lambda_{M^{\infty}(T_{[0,t-\epsilon)}(\Phi),C)}}(\epsilon,v)-1}{\epsilon}\right]\\
&\leq-\textbf{E}\Biggl[\prod_{z \in M^{\infty}(T_{[0,t-\epsilon)}(\Phi),C)}v(z)
\Biggl(\int_{\mathbb{R}^{2}}\textbf{1}_{v(x)> 1}(1-v(x))\prod_{z \in M^{\infty}(T_{[0,t-\epsilon)}(\Phi),C)}\\
&(1-h(x,z))\Lambda(dx)+\int_{\mathbb{R}^{2}}\textbf{1}_{v(x)\leq1}(1-v(x))\prod_{z \in M^{\infty}(T_{[0,t-\epsilon)}(\Phi),C)}(1-h(x,z))\\
&e^{-\epsilon\int_{\mathbb{R}^{2}}(1-h(x,y))\prod_{z \in M^{\infty}(T_{[0,t-\epsilon)}(\Phi),C)}(1-h(y,z))\Lambda(dy)}\Lambda(dx)+\\
& -\vartheta\Biggl(\epsilon,\int_{\mathbb{R}^{2}}(1-v(x))\prod_{z \in M^{\infty}(T_{[0,t-\epsilon)}(\Phi),C)}(1-h(x,z))\Lambda(dx)\Biggr)\Biggr)\Biggr].
\end{align*}
Now, we can let $\epsilon$ go to $0$ and use the continuity of $f_{\infty,\Lambda_{M^{\infty}(T_{[0,t-\epsilon)}(\Phi),C)}}$ to get:
\begin{align*}
 \lim_{\epsilon\rightarrow 0}\frac{f(t,v)-f(t-\epsilon,v)}{\epsilon}=-\int_{\mathbb{R}^{2}}f_{\infty,\Lambda}(t,H(v,x))(1-v(x))\Lambda(dx).
\end{align*}
 
\subsection{Differential Equation for $k$-Mat\'ern}\label{SS: k-Matern eq}
For general $k$, in order to get the differential equation of $f_{k,\Lambda}(t,v)$,
we need some observations on the structure of $M^{k}(\Xi,C)$ for an abitrary p.p. $\Xi$.
To each point in $\Xi$, we associate an infinite binary sequence $\textbf{e}(x)= e_{1}(x)e_{2}(x)$ $\cdots e_{i}(x)\cdots$,
where $e_{i}(x)=\textbf{1}_{x \in M^{i}(\Xi,C)}$ as defined in Section \ref{S:model}.
From the monotonicity property of $M^{i}$ and Proposition \ref{P: k Matern convergence to CSMA},
we observe that:
\begin{itemize}
 \item $\{e_{2i}(x)\}_{i \in \mathbb{N}}$ is a decreasing sequence;
 \item $\{e_{2i+1}(x)\}_{i \in \mathbb{N}}$ is an increasing sequence;
 \item $\lim\limits_{i \rightarrow \infty}e_{2i}(x)=\lim\limits_{i \rightarrow \infty}e_{2i+1}(x)= e_{\infty}(x)=\textbf{1}_{x \in M^{\infty}(\Xi,C)}$.
\end{itemize}
In order to investigate the joint structure of the $M^{i}(\Xi,C)$ processes with $i \leq k$,
we extract the first $k$ bits of $\textbf{e}(x)$ for each point $x$ in $\Xi$, 
namely $\textbf{e}^{k}(x)=e_{1}(x)\cdots e_{k}(x)$. From the above three observations,
$\textbf{e}^{k}(x)$ can only take one of the $k+1$ values $E^{k}_{j}, j\in \llbracket0,k\rrbracket$, where $E^{k}_{j}$ is defined as:
\begin{align*}
 E^{k}_{j}=\left\{\begin{array}{l r }
  \{01\}^{j}0^{k-2j}& \qquad j\in\llbracket0,\lfloor\frac{k}{2}\rfloor\rrbracket\\
  \{01\}^{k-j}1^{2j-k}& \qquad j\in\llbracket\lfloor\frac{k}{2}\rfloor+1,k\rrbracket
  \end{array}\right..
\end{align*}
 Examples of $E^{k}_{j}$ in the cases $k=5$ and $k=6$ are:
\begin{align*}
 \begin{array}{c c c c}
  E^{5}_{0}=00000 &E^{5}_{3}=01011 &E^{6}_{0}=000000 &E^{6}_{4}=010111\\
  E^{5}_{1}=01000 &E^{5}_{4}=01111 &E^{6}_{1}=010000 &E^{6}_{5}=011111\\
  E^{5}_{2}=01010 &E^{5}_{5}=11111 &E^{6}_{2}=010100 &E^{6}_{6}=111111\\
		  &	           &E^{6}_{3}=010101
\end{array}
\end{align*}
The properties of $E^{k}_{j}$ are:
\begin{itemize}
 \item $|E^{k}_{j}|_{1}=j$, where $|e|_{1}$ is the number of bits equal to $1$  in $e$.
 \item $\{E^{k}_{j}\}_{j=0}^{k}$ is an increasing sequence for the natural order of binary sequences.
 \item Let $\mathbb{OR}$ be the bit-by-bit ``or'' operator, then $\mathbb{OR}_{j \in J}E^{k}_{j}= E^{k}_{\max_{j \in J}j}$.
\end{itemize}
We can then define $M^{k}_{j}(\Xi,C)= \{x \in \Xi \textnormal{ s.t. } e^{k}(x)=E^{k}_{j}\}$ and get:
\begin{itemize}
  \item For each $k$, $M^{k}_{j}(\Xi,C), j\in \llbracket \lfloor\frac{k+1}{2}\rfloor,k\rrbracket$ forms a partition of $M^{k}(\Xi,C)$.
 \item For each $k$, $M^{k}_{j}(\Xi,C),j\in\llbracket0,k\rrbracket$ forms a partition of $\Xi$.
 \item For each $k$, $M^{k}_{j}(\Xi,C)=M^{k+1}_{j}(\Xi,C)$ for $j<\lfloor\frac{k}{2}\rfloor$. Since for $j<k/2$, $E^{k}_{j}$ is a prefix of $E^{k+1}_{j}$, we have  $M^{k}_{j}(\Xi,C) \supseteq M^{k+1}_{j}(\Xi,C)$. However, for $j<\lfloor\frac{k}{2}\rfloor$ and each $x$ in $M^{k}_{j}(\Xi,C)$, we have $e_{k}(x)=e_{k-1}(x)=0$. Thus,  $e_{k+1}(x)=0$. Hence, $E^{k+1}_{j}$ is the only possible extension of $E^{k}_{j}$ and then $M^{k}_{j}(\Xi,C)=M^{k+1}_{j}(\Xi,C)$ for $j<\lfloor\frac{k}{2}\rfloor$.
 \item For each $k$, $M^{k}_{j}(\Xi,C)=M^{k+1}_{j+1}(\Xi,C)$ for $j>\lfloor\frac{k}{2}\rfloor$. The reasoning is of the same type as above. The only difference is that $e_{k-1}(x)=e_{k}(x)=e_{k+1}(x)=1$.
 \item $M^{k}_{\lfloor\frac{k}{2}\rfloor}(\Xi,C)=M^{k+1}_{\lfloor\frac{k}{2}\rfloor}(\Xi,C)\cup M^{k+1}_{\lfloor\frac{k+2}{2}\rfloor}(\Xi,C)$. This is a consequence of the three above observations.\\
\end{itemize}
We now come back to our infinite timer marked $Poisson$ p.p. $\Phi$. The joint distribution of its $i$-Mat\'ern ($i=\llbracket0,k\rrbracket$) can be obtained from the functional 
\[g_{k,\Lambda}(t,v_{0},v_{1},\cdots,v_{k})= \textbf{E}\left[\prod_{j=0}^{k}\prod_{x \in M^{k}_{i}(T_{[0,t)}(\Phi),C)}v_{i}(x)\right],\]
since the joint generating functional of $M^{i}(T_{[0,t)}(\Phi),C)$, $i=\llbracket0,k\rrbracket$ can be expressed as:
\[\textbf{E}\left[\prod_{j=0}^{k}\prod_{x \in M^{i}(T_{[0,t)}(\Phi),C)}v_{i}(x)\right]= g_{k,\Lambda}(t,u_{0},u_{1},\cdots,u_{k}),\]
where $u_{i}(x)= \prod_{j=0}^{k} v_{j}(x)^{E^{k}_{i}(j)}$ and $E^{k}_{i}(0)E^{k}_{i}(1)\cdots E^{k}_{i}(k)$ is the binary representation of $E^{k}_{i}$. In the above definition, $v_{i}$ are functions from $\mathbb{R}^{2}$ to $[0,1]$ satisfying:
 \begin{align*}
  \int_{\mathbb{R}^{2}}|1-v_{i}(x)|\Lambda(dx)<\infty.
 \end{align*}

In particular, the generating functional of the $k$-Mat\'ern model, $f_{k,\Lambda}(t,v)$, can be obtained by taking $v_{0}=v_{1}=\cdots=v_{\lfloor\frac{k-1}{2}\rfloor}=\textbf{1}$ (the constant function which takes value $1$) and $v_{\lfloor\frac{k+1}{2}\rfloor}=\cdots=v_{k}=v$:
\[f_{k,\Lambda}(t,v)=g_{k,\Lambda}(t,\textbf{1},\textbf{1},\cdots,\textbf{1},v,\cdots,v) \qquad \textnormal{ there are $\lfloor\frac{k+1}{2}\rfloor$ $\textbf{1}$'s}.\] 

We now show that:
\begin{theorem}\label{T:k-Matern eq}
 The functional $g_{k,\Lambda}$ satisfies the following equation:
\begin{align}
  \nonumber&\frac{dg_{k,\Lambda}}{dt}(0,v_{0},\cdots,v_{k})=1\\
  \nonumber&\frac{dg_{k,\Lambda}}{dt}(t,v_{0},\cdots,v_{k})=-\int\limits_{\mathbb{R}^{2}}(1-v_{0}(x))g_{k,\Lambda}(t,v_{0},\cdots,v_{k})\Lambda(dx)\\
  \nonumber&-\sum_{i=0}^{\lfloor\frac{k}{2}\rfloor-1}\int\limits_{\mathbb{R}^{2}}(v_{i}(x)-v_{i+1}(x))g_{k,\Lambda}(t,v_{0,\cdots,H(v_{k-i},x),\cdots,H(v_{k},x)})\Lambda(dx)\\
  &-\sum_{i=\lfloor\frac{k}{2}\rfloor}^{k-1}\int\limits_{\mathbb{R}^{2}}(v_{i}(x)-v_{i+1}(x))g_{k,\Lambda}(t,v_{0,\cdots,H(v_{k-i-1},x),\cdots,H(v_{k},x)})\Lambda(dx).\label{E:k-Matern eq}
\end{align}
\end{theorem}
\emph{Proof} See appendix \ref{A:Proof T k Matern}\\
Now we can express the derivative of the functional $f_{k,\Lambda}$ as follows:
\begin{corollary}
 \begin{align*} 
&\frac{df_{2k,\Lambda}(t,v)}{dt}= -\int_{\mathbb{R}^{2}}g_{2k,\Lambda}(\textbf{1},\cdots,\textbf{1},v,H(v,x),\cdots,H(v,x))(1-v(x))\Lambda(dx)\\ 
&\qquad \textnormal{with $k$ $  \textbf{1}$'s} \\
&\frac{df_{2k+1,\Lambda}(t,v)}{dt}= -\int_{\mathbb{R}^{2}}g_{2k+1,\Lambda}(\textbf{1},\cdots,\textbf{1},H(\textbf{1},x),H(v,x),\cdots,H(v,x))(1-v(x))\Lambda(dx)\\ 
&\qquad \textnormal{with $k$ $  \textbf{1}$'s}.
 \end{align*}
\end{corollary}
\emph{Proof.}\\
It is sufficient to subtitute the formulas for $f_{2k,\Lambda}$ and $f_{2k+1,\Lambda}$ in Equation (\ref{E:k-Matern eq}).
\begin{flushright}
 $\Box$.
\end{flushright}

\section{Generating Functional of $\infty$-Mat\'ern under Palm Distribution} \label{S:Palm-Frechet}
The subject of this section is to study the generating functional of the $\infty$-Mat\'ern model under its Palm distribution. For any p.p., we show that with a proper parameterization of the function $v$, the generating functional under the Palm distribution can be expressed as the ``differential`` of the non-Palm version. Using this result and Theorem \ref{T:CSMA eq}, one can derive the differential equation that governs the evolution of the Palm version of the generating functional of the $\infty$-Mat\'ern model. We begin by recalling the formal definition of the Palm distribution.
\subsection{Radon-Nikodym Definition of the Palm Distribution}\label{SS:Palm}
When studying a p.p., it is sometime desirable to consider a point $x$ and condition on the fact that this p.p. contains $x$. The distribution of the p.p.  conditioned on this is called the Palm distribution w.r.t. $x$. In this subsection, we give the definition of Palm distributions for a general p.p.s. in terms of the Radon-Nikodym derivative of the Campbell measure \cite{Stoyan,DarVer1}. For this, we need to recall the formal definition of point proccesses.\\
Let $\mathbb{N}$ be  the set of sequences $\phi$ of points which are:
\begin{itemize}
 \item \emph{Locally finite}: for any bounded Borel set $B$ of $\mathbb{R}^{2}$, the number of points of $\phi$ in $B$ is finite.
 \item \emph{Simple}: each point appears only once in $\phi$.
\end{itemize}
Let $\mathcal{N}$ be the smallest sigma algebra that makes the functions $\phi \rightarrow \phi(B)$ measurable for all Borel sets $B$. A point process $\Xi$ is a measurable application from a probability space $(\Omega,\mathcal{A},\mathbb{P})$ to $(\mathbb{N},\mathcal{N})$.\\
For p.p. $\Xi$, let $\Lambda_{\Xi}$ be its intensity measure (i.e. $\Lambda_{\Xi}(B)=\textbf{E}[\Xi(B)]$ for each Borel set $B$). The Campbell measure is defined as:
\begin{align*}
 \mathcal{D}_{\Xi}(B,Y)=\textbf{E}[\Xi(B)\textbf{1}_{Y}],
\end{align*}
for all Borel sets $B$ and measurable sets $Y$ in $\mathcal{N}$. For each $Y$, it is easy to see that $\mathcal{D}_{\Xi}(.,Y)$ is absolutely continuous w.r.t. $\Lambda_{\Xi}(.)$ and hence admits a Radon-Nikodym density w.r.t. the later. Let us call this density $\textbf{P}_{x,\Xi}(Y)$, which satisfies:
\begin{align*}
 \mathcal{D}_{\Xi}(B,Y)= \int_{B}\textbf{P}_{x,\Xi}(Y)\Lambda_{\Xi}(dx) \qquad \textnormal{for any Borel set } B.
\end{align*}

It is not difficult to see that $\textbf{P}_{x,\Xi}(.)$ is actually a measure and that $\textbf{P}_{x,\Xi}(X)=1$ (see \cite{Stoyan}). Hence, $\textbf{P}_{x,\Xi}(.)$ is a probability measure and is called the Palm distribution of $\Xi$ w.r.t. $x$. In the special case where $\Xi$ is stationary, $\Lambda_{\Xi}= \lambda_{\Xi}\emph{l}$ ($\emph{l}$ is the Lesbegue measure) and:
\[\mathcal{D}_{\Xi}(B,Y)=\lambda_{\Xi}\int_{B}\textbf{P}_{x,\Xi}(Y)dx=\lambda_{\Xi}\int_{B}\textbf{P}_{o,\Xi}(Y_{-x})dx,\]
with $Y_{-x}:= \{\phi \textnormal{ s.t. } \phi+x \in Y\}$.\\
The following refined Campbell formula holds for all measurable functions $\gamma(x,\phi)$:
\begin{align}
& \nonumber \textbf{E}\left[\sum_{x \in \Xi} \gamma(x,\Xi)\right]=\int \gamma(x,\phi) \mathcal{C}_{\Xi}(dx,d\phi)=\int_{\mathbb{R}^{2}}\int \gamma(x,\phi) \textbf{P}_{x}(d\phi)\Lambda_{\Xi}(dx).\label{e:refined Campbell} 
\end{align}
  In the stationary case, we have:
\begin{align}
\textbf{E}\left[\sum_{x \in \Xi} \gamma(x,\Xi)\right] =\lambda\int_{\mathbb{R}^{2}}\int \gamma(x,\phi_{-x}) \textbf{P}_{o,\Phi}(d\phi)dx. 
\end{align}
More generally, we can define the $n$ fold Palm distribution of a p.p. $\Xi$ w.r.t. $x_{1},\cdots,x_{n}$ as the probability distribution of this p.p. conditioned on the fact that it contains the $n$ points $x_{1},\cdots,x_{n}$. First we need to define the $n$ fold Campbell measure as follows:
\begin{align*}
& \mathcal{C}^{n!}_{\Xi}(B_{1},\cdots,B_{n},Y)\\
& = \textbf{E}\left[\sum_{x_{1} \in \Xi\cap B_{1},\cdots,x_{n} \in \Xi \cap B_{n}; x_{1},\cdots,x_{n} \textnormal{ mutually different }}\textbf{1}\{\Xi \in Y\}\right].
\end{align*}
Let $\Lambda^{n!}_{\Xi}$ be the $n^{th}$ factorial moment measure of $\Xi$ \cite{Stoyan}. We can see that for each $Y$, the measure $\mathcal{C}^{n!}_{\Xi}$ is absolutely continuous w.r.t $\Lambda^{n!}_{\Xi}$. Hence it admits a Radon-Nikodym derivative w.r.t. $\Lambda^{n!}_{\Xi}$:
\begin{align*}
 \mathcal{C}^{n!}_{\Xi}(B_{1},\cdots,B_{n},Y)= \int_{B_{1}}\cdots\int_{B_{n}}\textbf{P}_{x_{1},\cdots,x_{n},\Xi}(Y)\Lambda^{n!}_{\Xi}(dx).
\end{align*}
One can then show that for each $n$-tuple $(x_{1},\cdots,x_{n})$, $\textbf{P}_{x_{1},\cdots,x_{n},\Xi}(.)$ is a probability measure which is called the $n$ fold Palm distribution of the p.p. $\Xi$ w.r.t. $x_{1},\cdots,x_{n}$. The following Campbell formula holds for all measurable functions $\gamma(x_{1},\cdots,x_{n},\phi)$:
\begin{align}
& \nonumber \textbf{E}[\sum_{x_{1},\cdots,x_{n} \in \Xi \textnormal{ mutually different}} \gamma(x_{1},\cdots,x_{n},\Xi)]\\
&= \nonumber \int \gamma(x_{1},..,x_{n},\phi) \mathcal{C}^{n!}_{\Xi}(dx_{1},\cdots,dx_{n},d\phi)\\
&= \int_{(\mathbb{R}^{2})^{n}}\gamma(x_{1},\cdots,x_{n},\phi)\textbf{P}_{x_{1},\cdots,x_{n},\Xi}(d\phi)\Lambda^{n!}_{\Xi}(dx_{1},..,dx_{n}).
\end{align}
\subsection{Generating Functionals under the Palm Distribution.}\label{SS:link}
Recall that for all p.p.s $\Xi$ and all functions $v$ satisfying (\ref{E: Gen def con}), the generating functional of $\Xi$ at $v$ is:
\begin{align*}
 G_{\Xi}(v)=\textbf{E}\left[\prod_{y\in \Xi}v(y)\right].
\end{align*}
The generating functional of $\Xi$ under its Palm distribution w.r.t. $x$ is defined as:
\begin{align*}
 G_{x,\Xi}(v)=\textbf{E}_{x,\Xi}\left[\prod_{y\in \Xi\setminus {x}}v(y)\right].
\end{align*}
Let $\textbf{x}^{n}=\{x_{1},\cdots,x_{n}\}$; the generating functionals of $\Xi$ under its $n$-fold Palm distribution w.r.t. $\textbf{x}^{n}$ is defined as:
\begin{align*}
 G_{\textbf{x}^{n},\Xi}(v)=\textbf{E}_{\textbf{x}^{n},\Xi}\left[\prod_{y\in \Xi\setminus \textbf{x}^{n}}v(y)\right].
\end{align*}
Note that for all pairs of functions $v$, $u$ such that $v$ and $1-u$ satisfy (\ref{E: Gen def con}), and for all $s>0$, we have that $u+sv$ satisfies:
\begin{align*}
 \int_{\mathbb{R}^{2}}|1-(v(x)+su(x))|\Lambda(dx) \leq \int_{\mathbb{R}^{2}}|1-v(x)|\Lambda(dx)+s\int_{\mathbb{R}^{2}}|u(x)|\Lambda(dx)< \infty.
\end{align*}
Thus the generating functional can be defined for $v+su$. The relation between $G_{\Xi}$ and $G_{\textbf{x}^{n},\Xi}$ is shown in the following proposition:
\begin{proposition}\label{P:Gen func Palm}
 For all pairs of functions $v$, $u$ such that $v$ and $1-u$ satisfy (\ref{E: Gen def con}), and for all $s>0$:
\begin{align*}
\frac{d^{n}}{d^{n}t}G_{\Xi}(v+tu)\Biggr|_{t=s}= \int_{(\mathbb{R}^{2})^{n}}u(x_{1})\cdots u(x_{n})G_{\textbf{x}^{n},\Xi}(v+su)\Lambda^{n!}_{\Xi}(dx_{1},\cdots,dx_{n}),
\end{align*}
given that the right-hand side is finite for all $n$.
\end{proposition}
\emph{Proof}
We have:
\begin{align*}
& G_{\Xi}(v+(t+\epsilon)u)\\
&= \textbf{E}\left[\prod_{x \in \Xi}(v(x)+(t+\epsilon)u(x))\right]\\
&= \textbf{E}\Bigl[\prod_{x \in \Xi}(v(x)+tu(x))+\sum_{n=1}^{\infty}\epsilon^{n}\sum_{(x_{1},\cdots,x_{n}) \in \Xi^{n!}}u(x_{1})\cdots u(x_{n})\\
&\prod_{x \in \Xi\setminus \{x_{1},\cdots,x_{n}\}}(v(x)+tu(x))\Bigr]\\
&= G_{\Xi}(v+tu)+\sum_{n=1}^{\infty}\epsilon^{n}\textbf{E}\Bigl[\hspace{-.2cm}\sum_{(x_{1},\cdots,x_{n})\in \Xi^{n!}}\hspace{-.2cm}u(x_{1})\cdots u(x_{n})\hspace{-.2cm}\prod_{x \in \Xi\setminus \{x_{1},\cdots,x_{n}\}}\hspace{-.2cm}(v(x)+tu(x))\Bigr]
\end{align*}
Consider $u(x_{1})\cdots u(x_{n})\prod_{x \in \Xi\setminus \{x_{1},\cdots,x_{n}\}}(v(x)+tu(x))$ as a function of the $n$ points $x_{1},\cdots,x_{n}$ and the p.p. $\Xi$; we can apply the $n$ fold Campbell formula to get:
\begin{align*}
 &\textbf{E}\Biggl[\sum_{x_{1},\cdots,x_{n} \textnormal{ pairwise different }\in \Xi}u(x_{1})\cdots u(x_{n})\prod_{x \in \Xi\setminus \{x_{1},\cdots,x_{n}\}}(v(x)+tu(x))\Biggr]\\
 &=\frac{1}{n!}\int_{(\mathbb{R}^{2})^{n}}u(x_{1})\cdots u(x_{n})G_{\textbf{x}^{n},\Xi}(v+tu)\Lambda^{n!}_{\Xi}(dx_{1},\cdots,dx_{n}).
\end{align*}
The last line is finite by assumption. 
So we have the following expansion of $G_{\Xi}(v+(t+\epsilon)u)$ around $t$:
\begin{align*}
& G_{\Xi}(v+(t+\epsilon)u)\\
&=G_{\Xi}(v+tu)+\sum_{n=1}^{\infty}\frac{\epsilon^{n}}{n!}\int_{(\mathbb{R}^{2})^{n}}u(x_{1})\cdots u(x_{n})G_{\textbf{x}^{n},\Xi}(v+tu)\Lambda^{n!}_{\Xi}(dx_{1},\cdots,dx_{n}).
\end{align*}
This expansion shows that $G_{\Xi}(v+tu)$ is infinitely differentiable in $t$. Then, by comparing with the Taylor expansion, we can conclude that:
\begin{align*}
\frac{d^{n}}{d^{n}t}G_{\Xi}(v+tu)\Biggr|_{t=s}= \int_{(\mathbb{R}^{2})^{n}}u(x_{1})\cdots u(x_{n})G_{\textbf{x}^{n},\Xi}(v+su)\Lambda^{n!}_{\Xi}(dx_{1},\cdots,dx_{n}).
\end{align*}
\begin{flushright}
 $\Box$.
\end{flushright}

\subsection{Palm Distribution of the $\infty$-Mat\'ern Model}\label{SS:CS Palm gen}
In this subsection, we determine the Palm distribution of the $\infty$-Mat\'ern model defined in Section \ref{S:model} by its generating functional. In particular, we use the results in Subsection \ref{SS:link} to derive the generating functional of this model under its Palm distribution. Let\\
\[ f^{y}_{\infty,\Lambda}(t,v)= \textbf{E}_{y,M^{\infty}(T_{[0,t)}(\Phi),C)}\left[\prod_{x \in M^{\infty}(T_{[0,t)}(\Phi),C)\setminus y }v(x)\right]\]
be the generating functional of the p.p. $M^{\infty}(T_{[0,t)}(\Phi),C)$ under its Palm distribution with respect to $y$.
\begin{proposition}\label{P:Palm gen infty Matern}
 For $\Lambda$-almost all $y$, for all functions $v$ satisfying (\ref{E: Gen def con}) and all $t>0$, the funtional $f^{y}_{\infty,\Lambda}(t,v)$ satisfies the functional equation:
\begin{align*}
& f^{y}_{\infty}(t,v)=\frac{\int_{0}^{t}f_{\infty,\Lambda}(\tau,H(v,y))d\tau}{m^{1}_{\infty,\Lambda}(t,y)}\\
&-\int_{0}^{t}\int_{\mathbb{R}^{2}}f^{y}_{\infty,\Lambda}(t,H(v,x))(1-v(x))(1-h(x,y))\frac{m^{1}_{\infty,\Lambda}(\tau,y)}{m^{1}_{\infty,\Lambda}(t,y)}\Lambda(dx)d\tau,
\end{align*}
where $H(.,.)$ is defined in (\ref{E:h transform}), $f_{\infty,\Lambda}$ is defined in (\ref{E: F def}) and
\begin{align*}
 m^{1}_{\infty,\Lambda}(t,x)= \int_{0}^{t}f_{\infty,\Lambda}(\tau,H(1,x))d\tau,
\end{align*}
is the Radon-Nikodym derivative of the first moment measure of the p.p.\\
 $M^{\infty}(T_{[0,t)}(\Phi),C)$ w.r.t. $\Lambda$.
\end{proposition}
\emph{Proof}
 As $M^{\infty}(T_{[0,t)}(\Phi),C)$ is dominated by $T_{[0,t)}(\Phi)$, we have that for all $v$, $u$ such that $v$ and $1-u$ satisfy the condition (\ref{E: Gen def con}),
\begin{align}
\nonumber& \textbf{G}_{\textbf{y}^{n},M^{\infty}(T_{[0,t)}(\Phi),C)}\left[v+su\right]\leq \textbf{G}_{\textbf{y}^{n},M^{\infty}(T_{[0,t)}(\Phi),C)}\left[v'+su\right]\\
&\leq \textbf{G}_{\textbf{y}^{n},T_{[0,t)}(\Phi)}\left[v'+su\right]=e^{t\int_{\mathbb{R}^{2}}|1-v(x)|\Lambda(dx)+ts\int_{\mathbb{R}^{2}}|u(x)|\Lambda(dx)} \label{E:Poiss bound}, 
\end{align}
where $v'(x)=1+|1-v(x)|$ for all $x$. Thus, we can easily verify that the p.p. $M^{\infty}(T_{[0,t)}(\Phi),C)$ satisfies the condition of Proposition \ref{P:Gen func Palm}.\\
Now let us begin with the computation of $m^{1}_{\infty,\Lambda}(t,x)$. From Proposition \ref{P:Gen func Palm}, we have,
\begin{align*}
& \frac{d}{ds}f_{\infty,\Lambda}(t,1+s\textbf{1}_{B})\Bigr|_{s=0}\\
&=\int_{B}f^{x}_{\infty,\Lambda}(t,1)m^{1}_{\infty,\Lambda}(t,x)\Lambda(dx)=\int_{B}m^{1}_{\infty,\Lambda}(t,x)\Lambda(dx)
\end{align*}
for all Borel sets $B$. Moreover, by Theorem \ref{T:CSMA eq}:
\begin{align*}
 & f_{\infty,\Lambda}(t,1+s\textbf{1}_{B})= 1+\int_{0}^{t}\int_{B}f_{\infty,\Lambda}(\tau,H(1,x)+sH(\textbf{1}_{B},x))s\Lambda(dx)d\tau.
\end{align*}
Hence:
\begin{align*}
 \int_{B}m^{1}_{\infty,\Lambda}(t,x)\Lambda(dx)=\frac{d}{ds}\int_{0}^{t}\int_{B}f_{\infty,\Lambda}(\tau,H(1,x)+sH(\textbf{1}_{B},x))s\Lambda(dx)d\tau\Bigl|_{s=0}.
\end{align*}
We want to interchange the derivative and the integration in the left hand side. A sufficient condition for this interchange is:
\begin{align}
 \int_{0}^{t}\int_{B}\frac{d}{ds}\left|f_{\infty,\Lambda}(\tau,H(1,x)+sH(\textbf{1}_{B},x))s\right|\Lambda(dx)d\tau <\infty \textnormal{ }\forall s>0 .\label{e:interchange cond}
\end{align}
Recall that: $\mathcal{N}= \sup_{x \in \mathbb{R}^{2}}\int_{\mathbb{R}^{2}}h(x,z)\Lambda(dz)$.
Using inequality (\ref{E:Poiss bound}) and Proposition \ref{P:Gen func Palm}, for all $s$, we have:
\begin{align*}
& \left|\frac{d}{ds}f_{\infty,\Lambda}(\tau,H(1,x)+sH(\textbf{1}_{B},x))\right|\\
&=\left|\int_{B}f^{y}_{\infty,\Lambda}(\tau,H(1,x)+.sH(\textbf{1}_{B},x))m^{1}_{\infty,\Lambda}(\tau,y)\Lambda(dy)\right|\\
&\leq \left|\int_{B}e^{\tau\int_{\mathbb{R}^{2}}h(x,z)\Lambda(dz)+\tau s\int_{B}(1-h(x,z))\Lambda(dz)}\tau\Lambda(dy)\right|\\
&\leq \left|\int_{B}e^{\tau(\mathcal{N}+s\Lambda(B))}\tau\Lambda(dx)\right|=\tau\Lambda(B)e^{\tau(\mathcal{N}+s\Lambda(B))}.
\end{align*}
In the third line we have used the fact that $m_{\infty,\Lambda}(\tau,x)<\tau$ for all $x$ (as the $\infty$-Mat\'ern model is a thinning of $\Phi$).
Hence:
\begin{align*}
&\left| \frac{d}{ds}\left(f_{\infty,\Lambda}(\tau,H(1,x)+sH(\textbf{1}_{B},x))s\right)\right|\\
&\leq\left|f_{\infty,\Lambda}(\tau,H(1,x)+sH(\textbf{1}_{B},x))\right|+s\left|\frac{d}{ds}f_{\infty,\Lambda}(\tau,H(1,x)+sH(\textbf{1}_{B},x))\right|\\
&\leq e^{\tau(\mathcal{N}+s\Lambda(B))}+s\tau\Lambda(B)e^{\tau(N+s\Lambda(B))}=(1+s\tau\Lambda(B))e^{\tau(\mathcal{N}+s\Lambda(B))}. 
\end{align*}
We have:
\begin{align*}
&\int_{0}^{t}\int_{B}(1+s\tau\Lambda(B))e^{\tau(\mathcal{N}+s\Lambda(B))}\Lambda(dx)d\tau \leq t\Lambda(B)(1+st\Lambda(B))e^{t(\mathcal{N}+s\Lambda(B)}<\infty,
\end{align*}
where we have used the fact that the function under the integral is an increasing function. This gives the condition (\ref{e:interchange cond}) directly.
So it is legitimate to interchange the differentiation and the integration to get:
\begin{align*}
 &\int_{B}m^{1}_{\infty,\Lambda}(t,x)\Lambda(dx)\\
 &=\int_{0}^{t}\int_{B}\frac{d}{ds}\left(f_{\infty,\Lambda}(\tau,H(1,x)+sH(\textbf{1}_{B},x))s\right)\Lambda(dx)d\tau\Bigl|_{s=0}\\
 &=\int_{0}^{t}\int_{B}f_{\infty,\Lambda}(\tau,H(1,x)+sH(\textbf{1}_{B},x))\Lambda(dx)d\tau\Bigl|_{s=0}+\\
&\int_{0}^{t}\int_{B}s\frac{d}{ds}f_{\infty,\Lambda}(\tau,H(1,x)+sH(\textbf{1}_{B},x))\Lambda(dx)d\tau\Bigl|_{s=0}\\
 &=\int_{0}^{t}\int_{B}f_{\infty,\Lambda}(\tau,H(1,x))\Lambda(dx)d\tau.
\end{align*}
As this equality is true for all Borel sets $B$, we must have:
\[m^{1}_{\infty,\Lambda}(t,x)=\int_{0}^{t}f_{\infty,\Lambda}(\tau,H(1,x))d\tau,\]
$\Lambda$-almost everywhere.\\
Now we can compute $f^{y}_{\infty,\Lambda}(t,v)$. By Proposition \ref{P:Gen func Palm}, we have:
\begin{align*}
 \frac{d}{ds}f_{\infty,\Lambda}(t,v+s\textbf{1}_{B})\Bigr|_{s=0}=\int_{B}f^{y}_{\infty,\Lambda}(t,v)m^{1}_{\infty,\Lambda}(t,y)\Lambda(dy).
\end{align*}
By Theorem \ref{T:CSMA eq}:
\begin{align*}
& f_{\infty,\Lambda}(t,v+s\textbf{1}_{B})\\
&=1-\int_{0}^{t}\int_{\mathbb{R}^{2}}f_{\infty,\Lambda}(\tau,H(v,x)+sH(\textbf{1}_{B},x))(1-v(x)-s\textbf{1}_{B}(x))\Lambda(dx)d\tau.
\end{align*}
By the same argument as in the computation of the intensity measure, we have:
\begin{align*}
 &\int_{B}f^{y}_{\infty,\Lambda}(t,v)m^{1}_{\infty,\Lambda}\Lambda(dx)\\
&=-\frac{d}{ds}\left(\int_{0}^{t}\int_{\mathbb{R}^{2}}f_{\infty,\Lambda}(\tau,H(v,x)+sH(\textbf{1}_{B},x))(1-v(x)-s\textbf{1}_{B}(x))\Lambda(dx)d\tau\right)\Biggr|_{s=0}\\
&=-\int_{0}^{t}\int_{\mathbb{R}^{2}}\frac{d}{ds}\left(f_{\infty,\Lambda}(\tau,H(v,x)+sH(\textbf{1}_{B},x))(1-v(x)-s\textbf{1}_{B}(x))\right)\Bigr|_{s=0}\Lambda(dx)d\tau.
\end{align*}
Moreover:
\begin{align*}
 &\frac{d}{ds}\left(f_{\infty,\Lambda}(\tau,H(v,x)+sH(\textbf{1}_{B},x))(1-v(x)-s\textbf{1}_{B}(x))\right)\Bigr|_{s=0}\\
 &=\frac{d}{ds}(f_{\infty,\Lambda}(\tau,H(v,x)+sH(\textbf{1}_{B},x))\Bigr|_{s=0}(1-v(x))-\textbf{1}_{B}(x)f_{\infty,\Lambda}(\tau,H(v,x)).
\end{align*}
By Proposition \ref{P:Gen func Palm}, for all $x$:
\begin{align*}
  &\frac{d}{ds}(f_{\infty,\Lambda}(\tau,H(v,x)+sH(\textbf{1}_{B},x))\Bigr|_{s=0}\\
&=\int_{B}f^{y}_{\infty,\Lambda}(\tau,H(v,x))(1-h(x,y)))m^{1}_{\infty,\Lambda}(\tau,y)\Lambda(dy).
\end{align*}
Hence:
\begin{align*}
 &\int_{B}f^{y}_{\infty,\Lambda}(t,v)m^{1}_{\infty,\Lambda}(t,y)\Lambda(dx)\\
&=-\int_{0}^{t}\int_{\mathbb{R}^{2}}\int_{B}f^{y}_{\infty,\Lambda}(\tau,H(v,x))(1-h(x,y))m^{1}_{\infty,\Lambda}(\tau,y)\Lambda(dy)\Lambda(dx)d\tau\\
&+\int_{0}^{t}\int_{B}f_{\infty,\Lambda}(\tau,H(v,x))\Lambda(dx)d\tau\\
&=-\int_{B}\int_{0}^{t}\int_{\mathbb{R}^{2}}f^{y}_{\infty,\Lambda}(\tau,H(v,x))(1-h(x,y)))\Lambda(dx) m^{1}_{\infty,\Lambda}(\tau,y)d\tau\Lambda(dy)\\
&+\int_{B}\int_{0}^{t}f_{\infty,\Lambda}(\tau,H(v,y))d\tau\Lambda(dy).
\end{align*}
As this is true for any Borel set $B$, we must have:
\begin{align*}
& f^{y}_{\infty}(t,v)=\frac{\int_{0}^{t}f_{\infty,\Lambda}(\tau,H(v,y))d\tau}{m^{1}_{\infty,\Lambda}(t,y)}\\
&-\int_{0}^{t}\int_{\mathbb{R}^{2}}f^{y}_{\infty,\Lambda}(\tau,H(v,x))(1-v(x))(1-h(x,y))\frac{m^{1}_{\infty,\Lambda}(\tau,y)}{m^{1}_{\infty,\Lambda}(t,y)}\Lambda(dx)d\tau
\end{align*}
$\Lambda$- almost everywhere.
\begin{flushright}
 $\Box$
\end{flushright}
By extending this framework, we can recursively compute the $n$-th factorial moment measures and the generating functional under the $n$ fold Palm distribution of the $\infty$-Mat\'ern model. In these developments, for brevity we adopt the following notation, which extends that used in Proposition \ref{P:Gen func Palm}:
\begin{itemize}
 \item $\textbf{s}^{n}$ should be read $\{s_{1},\cdots,s_{n}\}$.
 \item $\textbf{s}^{n}_{-i}$ should be read $\{s_{1},\cdots,s_{i-1},s_{i+1},\cdots,s_{n}\}$.
 \item $\prod\textbf{s}^{n}=\prod_{i=1}^{n}s_{i}$ and $\prod\textbf{s}_{-i}^{n}=\prod_{j=1}^{i-1}s_{j}\prod_{j=i+1}^{n}s_{j}$.
 \item $\int_{\textbf{B}^{n}}=\int_{B_{1}}\cdots\int_{B_{n}}$ and $\int_{\textbf{B}_{-i}^{n}}=\int_{B_{1}}\cdots\int_{B_{i-1}}\int_{B_{i+1}}\cdots\int_{B_{n}}$
 \item $\Lambda^{n}(d\textbf{x}^{n})=\prod_{j=1}^{n}\Lambda(dx_{j})$ and $\Lambda^{n-1}(d\textbf{x}^{n}_{-i})=\prod_{j=1}^{i-1}\Lambda(dx_{j})\prod_{j=i+1}^{n}\Lambda(dx_{j})$.
 \item $\textbf{s}^{n}=t$ means $s_{i}=t$ for $i \in \llbracket 1,n \rrbracket$
\end{itemize}
Let 
\begin{align*}
 f^{\textbf{y}^{n}}_{\infty,\Lambda}(t,v)= \textbf{E}_{y_{1},\cdots,y_{n}}\left[\prod_{x \in M^{\infty}(T_{[0,t)}(\Phi),C)\setminus \{y_{1},\cdots,y_{n}\}}\right],
\end{align*}
be the generating functional of the $\infty$-Mat\'ern model under its $n$ fold Palm distribution w.r.t. $\{y_{1},\cdots,y_{n}\}$. We also let $m^{n}_{\Lambda}(t,y_{1},\cdots,y_{n})$ be the $n$-factorial density of $M^{\infty}(T_{[0,t)}(\Phi),C)$ (i.e. the Radon-Nikodym derivative of the $n$-factorial moment measure of $M^{\infty}(T_{[0,t)}(\Phi),C)$ w.r.t. the measure $\Lambda^{n}$).
Now we are in a position to state the main result:
\begin{proposition}\label{P:Palm gen Matern general}
 For $\Lambda^{n}$-almost all $n$-tuple $\textbf{y}^{n}=\{y_{1},\cdots,y_{n}\}$ of pairwise different points, for all functions $v$ satisfying (\ref{E: Gen def con}) and all $t>0$, the funtional $f^{\textbf{y}^{n}}_{\infty,\Lambda}(t,v)$ satisfies:
\begin{align*}
 &f^{\textbf{y}^{n}}_{\infty}(t,v)=\sum_{i=1}^{n}\int_{0}^{t}f^{\textbf{y}^{n}_{-j}}_{\infty,\Lambda}(\tau,H(v,y_{j}))\prod_{i \in\llbracket 1,n\rrbracket,i \neq j}(1-h(y_{i},y_{j}))\frac{m^{n-1}_{\infty,\Lambda}(\tau,\textbf{y}^{n}_{-j})}{m^{n}_{\infty,\Lambda}(t,\textbf{y}^{n})}d\tau\\
&-\int_{0}^{t}\int_{\mathbb{R}^{2}}f^{\textbf{y}^{n}}_{\infty,\Lambda}(\tau,H(v,x))\prod_{i=1}^{n}(1-h(y_{i},x))(1-v(x))\frac{m^{n}_{\infty,\Lambda}(\tau,\textbf{y}^{n})}{m^{n}_{\infty,\Lambda}(t,\textbf{y}^{n})}\Lambda(dx)d\tau.
\end{align*}
where:
\begin{align*}
 m^{n}_{\infty,\Lambda}(t,\textbf{x}^{n})= \sum_{j=1}^{n}\int_{0}^{t}f^{\textbf{x}^{n}_{-j}}_{\infty,\Lambda}(\tau,H(1,x_{i}))\prod_{i \in\llbracket 1,n\rrbracket,i \neq j}(1-h(x_{i},x_{j}))m^{n-1}_{\infty,\Lambda}(\tau,\textbf{x}^{n}_{-i})d\tau,
\end{align*}
is the $n$-factorial moment density of the p.p. $M^{\infty}(T_{[0,t)}(\Phi),C)$.
\end{proposition}
\emph{Proof.} See appendix \ref{A:Proof P Palm general}.\\

\section{Solution of the Differential Equation}\label{S:diff solution}
In Section \ref{S:Anal} we have shown that the functionals $f_{\infty,\Lambda}$ and $g_{k,\Lambda}$ can be characterized by the system of differential equations (\ref{E:CSMA eq}) and (\ref{E:k-Matern eq}) respectively. Now we will prove the converse: given a system of differential equations of the form (\ref{E:CSMA eq}) and (\ref{E:k-Matern eq}), under some mild conditions, each system of equations has one and only one positive solution. These solutions are $f_{\infty,\Lambda}$ and $g_{k,\Lambda}$ respectively. Here, we have to restrict ourselves to the case where $v$ is smaller than $1$. This result is still strong enough for many applications, as the distribution of a p.p. is uniquely determined by its generating functionals at functions $v$ smaller than $1$. 
Note that we are working with generating functionals, so they must be a priori positive. With this condition, we have the following proposition.
\begin{proposition}\label{P: CSMA eq unique solution}
The system of differential equations (\ref{E:CSMA eq}) has a unique positive solution.
\end{proposition}
\emph{Proof}.\\
We show that if a functional $f$ satisfies this equation,
then it can only have one value for each $t$ and $v$. To simplify the notation, we consider the case where $\Lambda$ is the Lebesgue measure. The argument is similar for the general case.
The proof will be by induction on $t$.
Fix a positive number $t_{0}<N^{-1}$. 
For $t\leq t_{0}$, we have:
\begin{equation}
\label{eq:int}
f(t,v)= 1-\int_{0}^{t}\int_{\mathbb{R}^{2}}f(\tau,H(v,x))(1-v(x))\lambda dxd\tau.
\end{equation}
By putting $H(v,y)$ in place of $v$ in the above inequatlity and by
using the fact that $f(t,v)\leq 1$ for all $t,v$, we get that:
\begin{align*}
1\geq f(t,H(v,y)) \geq 1-t\int_{\mathbb{R}^{2}}(1-H(v,y)(x))dx \qquad \forall y \in \mathbb{R}^{2}.
\end{align*}
When using this in (\ref{eq:int}), we get
\begin{eqnarray*}
&& \hspace{-.8cm}1-t\int\limits_{\mathbb{R}^{2}}(1-v(x))dx+ \frac{t^{2}}{2}\hspace{-.2cm}\int\limits_{(\mathbb{R}^{2})^{2}}(1-v(x))(1-H(v,x)(y))dxdy\\
&& \geq f(t,v)\geq 1-t\int\limits_{\mathbb{R}^{2}}(1-v(x))dx.
\end{eqnarray*}
By repeating this we have:
\begin{align*}
&1+\sum_{k=1}^{2n}(-1)^{k}\frac{t^{k}}{k!}\int\limits_{(\mathbb{R}^{2})^{k}}\prod_{i=1}^{k}(1-v(x_{i})\prod_{j=1}^{i-1}(1-h(x_{j},x_{i})))dx_{1}\cdots dx_{k} \geq f(t,v) \geq\\
&1+\sum_{k=1}^{2n-1}(-1)^{k}\frac{t^{k}}{k!}\int\limits_{(\mathbb{R}^{2})^{k}}\prod_{i=1}^{k}(1-v(x_{i})\prod_{j=1}^{i-1}(1-h(x_{j},x_{i})))dx_{1}\cdots dx_{k}.
\end{align*}
All we need now is to prove that the series:
\[1+\sum_{k=1}^{\infty}(-1)^{k}\frac{t^{k}}{k!}\int_{(\mathbb{R}^{2})^{k}}\prod_{i=1}^{k}(1-v(x_{i})\prod_{j=1}^{i-1}(1-h(x_{j},x_{i})))dx_{1}\cdots dx_{k}\]
 converges. First we see that the following inequalities hold for all $k$:
 \begin{align*}
  &\int_{(\mathbb{R}^{2})^{k}}\prod_{i=1}^{k}(1-v(x_{i})\prod_{j=1}^{i-1}(1-h(x_{j},x_{i})))dx_{1}\cdots dx_{k}\\
  &\leq	\int_{(\mathbb{R}^{2})^{k}}\prod_{i=1}^{k}(1-v(x_{i})+\sum_{j=1}^{i-1}v(x_{i})h(x_{j},x_{i})))dx_{1}\cdots dx_{k}\\
  &\leq \int_{(\mathbb{R}^{2})^{k}}\prod_{i=1}^{k}(1-v(x_{i})+\sum_{j=1}^{i-1}h(x_{j},x_{i})))dx_{1}\cdots dx_{k}\\
  &= \prod_{i=1}^{k}\left(\int_{\mathbb{R}^{2}}(1-v(x))dx+(i-1)N\right).
 \end{align*}
Thus, if $t\leq t_{0}<N^{-1}$ then 
  \begin{align*}
   &\lim_{k\rightarrow \infty}\left(\frac{t^{k}}{k!}\int\limits_{(\mathbb{R}^{2})^{k}}
\prod_{i=1}^{k}(1-v(x_{i})\prod_{j=1}^{i-1}(1-h(x_{j},x_{i})))dx_{1}\cdots dx_{k}\right)^{1/k}\\
   &\leq \lim_{k\rightarrow \infty}\left(\frac{t^{k}}{k!}\prod_{i=1}^{k}\left(
\int\limits_{\mathbb{R}^{2}}((1-v(x))dx+(i-1)N)\right)\right)^{1/k}= tN < 1.
  \end{align*}
So our series indeed converges and the value of $f(t,v)$ can be uniquely determined as:
\begin{eqnarray}
 f(t,v) = 1+\sum_{k=1}^{\infty}\frac{(-t)^{k}}{k!}\int\limits_{(\mathbb{R}^{2})^{k}}
\prod_{i=1}^{k}(1-v(x_{i})\prod_{j=1}^{i-1}(1-h(x_{j},x_{i}))dx_{1}\cdots dx_{k}. \label{eq:toutestla}
\nonumber
\end{eqnarray}
Now suppose that the value of the functional $f$ is determined up to $t\leq mt_{0}$. 
By the very same technique we can prove that the series
\begin{align*}
 &  f(mt_{0},v)+\sum_{k=1}^{\infty}(-1)^{k}\frac{(t-mt_{0})^{k}}{k!}\int_{(\mathbb{R}^{2})^{k}}\prod_{i=1}^{k}f\left(mt_{0},v(.)\prod_{j=1}^{i}(1-h(x_{j},.))\right)\\ 
 & \left(1-v(x_{i})\prod_{j=1}^{i-1}(1-h(x_{j},x_{i}))\right)dx_{1}\cdots dx_{k}
\end{align*}
converges and we can determine $f(t,v)$ as:
\begin{align*}
 & f(t,v) =  f(mt_{0},v) \\
&+\sum_{k=1}^{\infty}(-1)^{k}\frac{(t-mt_{0})^{k}}{k!}
\int\limits_{(\mathbb{R}^{2})^{k}}\prod_{i=1}^{k}f\left(mt_{0},v(.)\prod_{j=1}^{i}(1-h(x_{j},.))\right)\\
& \hspace{1cm}\left(1-v(x_{i})\prod_{j=1}^{i-1}(1-h(x_{j},x_{i}))\right)dx_{1}\cdots dx_{k}.
\end{align*}
This terminates our proof.
\begin{flushright}
 $\Box$
\end{flushright}
For the $k-$Mat\'ern model, due to the complexity of the differential equation involving $g_{k,\Lambda}$, we limit ourselves to a less quantitative result. Let us first notice that the series representation of $f_{\infty,\Lambda}$ is its Taylor expansion at $n\delta$ with convergence radius larger than $\delta$. This gives rise to the following observation:
\begin{lemma}\label{L: Talyor unique cond}
 If a function $f$ is infinitely differentiable and for all $t$, the convergence radius for Taylor expansion at $t$ of $f$ is larger than $\delta$ for some $\delta>0$, then, $f$ is uniquely determined.
\end{lemma}
\emph{Proof.}\\
This proof follows exactly the method used in the proof of Proposition \ref{P: CSMA eq unique solution}.
\begin{flushright}
 $\Box$
\end{flushright}
We then try to bound the $l^{th}$ derivative w.r.t. $t$ of $g_{k,\Lambda}(t,v_{0},..,v_{k})$.
\begin{lemma}\label{L: Mk diff bound}
 For all integers $k,l$ and all functions $v_{0},\cdots,v_{k}$ satisfying:
\begin{align*}
 \int_{\mathbb{R}^{2}}|1-v_{i}(x)|\Lambda(dx) < \infty \qquad \textnormal{for } i=\overline{0,k},
\end{align*}
we have:
\begin{align*}
\left|\frac{d^{l}}{d^{l}t}g_{k,\Lambda}(t,v_{0},\cdots,v_{k})\right| \leq \prod_{i=0}^{k-1}(M(v_{0},\cdots,v_{k})+iN),
\end{align*}
where:
\begin{align*}
 M(v_{0},\cdots,v_{k})&= \int_{\mathbb{R}^{2}}|1-v_{0}(x)|\Lambda(dx)+\sum_{i=0}^{k-1}\int_{\mathbb{R}^{2}}|v_{i}(x)-v_{i+1}(x)|\Lambda(dx) < \infty\\
 \mathcal{N} &=\sup_{y \in \mathbb{R}^{2}}\int_{\mathbb{R}^{2}}h(x,y)\Lambda(dx) < \infty.
\end{align*}
\end{lemma}
\emph{Proof.}\\
The proof is by induction on $l$. For $l=1$, we have from Theorem \ref{T:k-Matern eq} that:
\begin{align*}
 &\left|\frac{d}{dt}g_{k,\Lambda}(t,v_{0},\cdots,v_{k})\right|=\left|-\int_{\mathbb{R}^{2}}(1-v_{0}(x))g_{k,\Lambda}(t,v_{0},\cdots,v_{k})\Lambda(dx)\right.\\
&-\sum_{i=0}^{\lfloor\frac{k}{2}\rfloor-1}\int_{\mathbb{R}^{2}}(v_{i}(x)-v_{i+1}(x))g_{k,\Lambda}(t,v_{0},\cdots,H(v_{k-i},x),\cdots,H(v_{k},x))\Lambda(dx)\\\
&\left.-\sum_{i=\lfloor\frac{k}{2}\rfloor}^{k-1}\int_{\mathbb{R}^{2}}(v_{i}(x)-v_{i+1}(x))g_{k,\Lambda}(t,v_{0},\cdots,H(v_{k-i-1},x),\cdots,H(v_{k},x))\Lambda(dx)\right|\\
&\leq M,
\end{align*}
where we have used the fact that $g_{k,\Lambda}$ is always smaller than $1$ in the last inequality.\\
Now assuming that the lemma is true for $l$, we have again from Theorem \ref{T:k-Matern eq}:
\begin{align}
 \nonumber&\left|\frac{d^{l+1}}{d^{l+1}t}g_{k,\Lambda}(t,v_{0},\cdots,v_{k})\right|=\left|-\int_{\mathbb{R}^{2}}(1-v_{0}(x))\frac{d^{l}}{d^{l}t}g_{k,\Lambda}(t,v_{0},\cdots,v_{k})\Lambda(dx)\right.\\
\nonumber &-\sum_{j=0}^{\lfloor\frac{k}{2}\rfloor-1}\int_{\mathbb{R}^{2}}(v_{i}(x)-v_{i+1}(x))\frac{d^{l}}{d^{l}t}g_{k,\Lambda}(t,v_{0},\cdots,H(v_{k-i},x),\cdots,H(v_{k},x))\Lambda(dx)\\\
&\left.-\sum_{j=\lfloor\frac{k}{2}\rfloor}^{k-1}\int_{\mathbb{R}^{2}}(v_{i}(x)-v_{i+1}(x))\frac{d^{l}}{d^{l}t}g_{k,\Lambda}(t,v_{0},\cdots,H(v_{k-i-1},x),\cdots,H(v_{k},x))\Lambda(dx)\right|. \label{E:ind step}
\end{align}
Adopting the convention that $v_{-1}=1$, we have for all $j$:
\begin{align}
& \nonumber \int_{\mathbb{R}^{2}}|H(x,v_{j})(y)-H(x,v_{j+1})(y)|\Lambda(dy) \\
& =\int_{\mathbb{R}^{2}}(1-h(x,y))|v_{j}(y)-v_{j+1}(y)|\Lambda(dy)\leq \int_{\mathbb{R}^{2}}|v_{j}(y)-v_{j+1}(y)|\Lambda(dy)\\
& \nonumber \int_{\mathbb{R}^{2}}|v_{j}(y)-H(x,v_{j+1})(y)|\Lambda(dy) \\
& =\nonumber \int_{\mathbb{R}^{2}}|(v_{j}(y)-v_{j+1}(y))+h(x,y)v_{i+1}(y)|\Lambda(dy) \\
& \leq \nonumber \int_{\mathbb{R}^{2}}(|(v_{j}(y)-v_{j+1}(y))|+|h(x,y)v_{i+1}(y)|)\Lambda(dy)\\
& \leq \int_{\mathbb{R}^{2}}|(v_{j}(y)-v_{j+1}(y))|\Lambda(dy) +\mathcal{N}.
\end{align}
Thus, for all $j>0$
\begin{align*}
 &\left|\frac{d^{l}}{d^{l}t}g_{k,\Lambda}(t,v_{0},\cdots,v_{j},H(x,v_{j+1}),\cdots,H(x,v_{k}))\right|\\
 &\leq\prod_{j=0}^{l-1}(M(v_{0},\cdots,v_{j},H(x,v_{j+1}),\cdots,H(x,v_{k}))+j\mathcal{N})\\
 &\leq\prod_{j=1}^{l}(M(v_{0},\cdots,v_{j},v_{j+1},\cdots,v_{k})+j\mathcal{N}).
\end{align*}
Now, plugging this into Equation (\ref{E:ind step}) gives us the desired result.
\begin{flushright}
 $\Box$
\end{flushright}
Combining the above two lemmas, we can prove the following:
\begin{proposition}\label{P:Mk eq unique}
The functional $g_{k,\Lambda}(t,v_{0},\cdots,v_{k})$ is uniquely determined. 
\end{proposition}
From Theorem \ref{T:k-Matern eq} we know that $g_{k,\Lambda}$ is infinitely differentiable with respect to $t$. Further, by using Lemma \ref{L: Mk diff bound}, we can show that the convergence radius of the Taylor expansion of $g_{k,\Lambda}$ at $t$ is larger than $\delta$ with $\delta< \mathcal{N}^{-1}$. Then using Lemma \ref{L: Talyor unique cond} gives the conclusion.
\begin{flushright}
 $\Box$.
\end{flushright}
\begin{corollary}
 The functional $f_{k,\Lambda}$ is uniquely determined.
\end{corollary}

\section{Bounds}
In this section, we concentrate on the numerical method which can be used to evaluate the functional
$f_{\infty,\Lambda}$. For simplicity, we consider only the stationary case, where $\Lambda(dx)=\lambda dx$ and the notation is reduced to $f_{\infty}$. As stated earlier, the series expansion given in the proof of
Proposition \ref{P: CSMA eq unique solution}
is the Taylor expansion of $f_{\infty}$ at $t_{0}$.
The convergence radius of this expansion is $N^{-1}$
(the definition of $N$ is in (\ref{E:CSMA well def cond})).
The expressions in the proof show that one can evaluate
$f_{\infty}(t,v)$ in a recursive way as follows:\\
Suppose that $mt_{0}<t\leq (m+1)t_{0}$ and $t_{0}<N^{-1}$.
\begin{itemize}
 \item Write $f_{\infty}(t,v)$ as a series as in (\ref{eq:toutestla}).
 \item Estimate $f_{\infty}(mt_{0},v)$.
 \item Sum up the series in step 1 to some $k$ that sastifies the error requirement.
\end{itemize}
Nevertheless, this direct approach requires unacceptable computing
cost due to multiple nested integrations of complex integrands.
Even computing a single term of order $k$ in the series would be non-polynomial in $k$.\\
For this reason, we propose here a bound on $f_{\infty}(t,v)$ which holds for all $t$ and all $v$ smaller than $1$.
As we will see in the section on applications, this
still provides a good enough estimate.
This bound is based mainly on Theorem \ref{T:CSMA eq} and 
on the simple observation that $M^{\infty}(\Phi,C)$ is a thinning of $\Phi$.
\begin{corollary}\label{C: f bound}
 The generating functional of $M^{\infty}(\Phi)$ sastisfies the following inequalities:
\begin{align*}
 &e^{-t\int_{\mathbb{R}^{2}}(1-v(x))\lambda dx}\leq f_{\infty}(t,v)\leq\\
 &1- \int_{0}^{t}\int_{\mathbb{R}^{2}}\frac{1-e^{-\tau\int_{\mathbb{R}^{2}}(1-(1-h(x,y))v(y))\lambda dy}}
{\int_{\mathbb{R}^{2}}(1-(1-h(x,y))v(y))\lambda dy}(1-v(x))\lambda dxd\tau.
\end{align*}
\end{corollary}
\emph{Proof.}
The left-hand side inequality comes directly from the fact that $M^{\infty}(\Phi)$ is a thinnning of $\Phi$.
The right-hand side is obtained by plugging the left-hand side inequality into the expression
of $\frac{df_{\infty,\Lambda}}{dt}$ in Theorem \ref{T:CSMA eq} and then by integrating 
with respect to $\tau$ from $0$ to $t$.
\begin{flushright}
 $\Box$.
\end{flushright}

\section{Conclusion}
 In this paper, the distributions of p.p.s induced by several random packing models are investigated. A key ingredient of this study is the existence of a time dimension and of a quasi-Poisson property in the sense of Lemma \ref{L:limit 0}. Using these properties, the generating functionals of the p.p.s under consideration and their Palm versions are derived as the solutions of some differential equations. Bounds and numerical methods to compute these generating functionals are derived.\\
 We highlight here two of the most prominent unanswered questions:
\begin{enumerate}
 \item \emph{Solution of the differential equation:} Here we have just proved that the differental equations (\ref{E:CSMA eq}) and (\ref{E:k-Matern eq}) have a unique solution. For some applications, it is desirable to compute these solutions explicitly. Thus there is a strong need for a numerical method to solve these equations with high accuracy and reasonable complexity (at least faster than MC methods).
 \item \emph{Random packing with depacking:} In the model we considered, whenever a point is packed, it stays there forever. Sometime it is also useful to consider de-packing, i.e. a packed point is removed from the packed p.p. This kind of model arises in many disciplines. For example one can consider the CSMA protocol for wireless network. Suppose that each transmitter has a random amount of data to transmit. The more data, the longer it takes to transmit. A transmitter departs from the network once it has transmitted all its data. The random packing model with depacking is hence a natural model to describe this kind of CSMA network.\\
In fact, the random packing model with de-packing can be considered as a spatial birth and death process \cite{GK06}. The study of this kind of model leads to two main difficulties: (i) it is not clear whether such a process is well-defined, (ii) we do not know whether this model has the quasi-Poisson property (in the sense of Lemma \ref{L:limit 0}).
\end{enumerate}

\bibliographystyle{ieeetr}
\bibliography{CSMA6}

\appendix
\section{ Proof of Theorem \ref{T:k-Matern eq}}\label{A:Proof T k Matern}
For brevity, we present here an informal proof which only captures the main ideas. Nevertheless, if one wants to formalize this proof, one can follow the same method as in Theorem \ref{T:CSMA eq}.\\
Let us consider a fixed real number $t$ and a ``thin layer'' $T_{[t,t+\epsilon)}(\Phi)$. This layer is so thin that there is almost no inner contention between points. A point $x$ in $T_{[t,t+\epsilon)}(\Phi)$ belongs to one of $k+1$ p.p.s $M^{k}_{i}(T_{[0,t+\epsilon)}(\Phi),C), i=\overline{0,k}$. Since there is no contention between points in $T_{[t,t+\epsilon)}(\Phi)$, the condition that $x$ belongs to $M^{k}_{i}(T_{[0,t+\epsilon)}(\Phi),C)$ depends only on the set $\mathcal{D}(x)$ of its contenders whose timer is smaller than $t$, as shown in the following claim:\\
  \emph{\textbf{Claim:}}\\
For $i\in \llbracket0,k-1\rrbracket$, the point $x$ belongs to $M^{k}_{i}(T_{[0,t+\epsilon)}(\Phi),C)$ iff \\
$\mathcal{D}(x) \cap M^{k-1}_{k-1-i}(T_{[0,t)}(\Phi),C)\neq \emptyset$ and $\mathcal{D}(x) \cap M^{k-1}_{j}(T_{[0,t)}(\Phi),C)= \emptyset$ for all $j>k-1-i$. Equivalently, $\textbf{e}^{k}(x)= E^{k}_{i}$ iff $\mathbb{OR}_{y \in \mathcal{D}(x)}\textbf{e}^{k-1}(y)= E^{k-1}_{k-1-i}$ .\\
  \emph{\textbf{Proof of Claim:}}\\
  Note that for all $j$, $\mathbb{OR}_{y \in \mathcal{D}(x)}e_{j}(y)=0$ iff $\mathcal{D}(x)\cap M^{j}(T_{[0,t+\epsilon)}(\Phi),C) = \emptyset$ , or equivalently iff $x \in M^{j+1}(T_{[0,t+\epsilon)}(\Phi),C)$. Thus for all $i\in \llbracket0,k\rrbracket$, we have:
\[\textbf{e}^{k}(x)= E^{k}_{i} \Longleftrightarrow \mathbb{OR}_{y \in \mathcal{D}(x)}\textbf{e}^{k-1}(y)= E^{k-1}_{k-1-i}.\]
Let $J= \{j \textnormal{ s.t. } M^{k-1}_{j}(T_{[0,t+\epsilon)}(\Phi),C) \cap \mathcal{D}(x) \neq \emptyset\}$. By the third property of $E^{k}_{j}$'s we have:
   \[E^{k-1}_{k-1-i}=\mathbb{OR}_{y \in \mathcal{D}(x)}e_{k-1}(y)=\mathbb{OR}_{j \in J}E^{k-1}_{j} = E^{k-1}_{\max_{j \in J}j}.\]
Then $\max_{j \in J}j = k-1-i$, and this justifies our claim.\\
Since $\mathcal{D}(x)$ is an independent thinning of $\Phi$ with thinning probability $h(x,y)$, the probability that $\mathcal{D}(x) \cap M^{k-1}_{j}(T_{[0,t)}(\Phi),C)= \emptyset$ conditioned on $T_{[0,t)}(\Phi)$ is:
 \[\mathcal{H}^{k-1}_{j}(x,t):=\prod_{y \in M^{k-1}_{j}(T_{[0,t)}(\Phi),C)}(1-h(x,y)).\]
Thus, conditioned on $T_{[0,t)}(\Phi)$ one can regard $T_{[t,\epsilon)}(M^{k}_{i}(T_{[0,t+\epsilon)}(\Phi),C))$, $i= \overline{0,k-1}$ as $k$ disjoint  independent thinnings of $T_{[t,t+\epsilon)}(\Phi)$ with the thinning probabilities:
 \[\left(1-\mathcal{H}^{k-1}_{k-1-i}(x,t)\right)\prod_{j=k-i}^{k-1}\mathcal{H}^{k-1}_{j}(x,t).\] 
So we have that:
\[T_{[t+\epsilon)}(M^{k}_{k}(T_{[0,t+\epsilon)}(\Phi),C))= T_{[t,t+\epsilon)}(\Phi)\setminus \left(\bigcup_{i=0}^{k-1}T_{[t,\epsilon)}(M^{k}_{i}(T_{[0,t+\epsilon)}(\Phi),C))\right),\]
 is an independent thinning of $T_{[t,t+\epsilon)}(\Phi)$ with thinning probability:
\[\prod_{j=0}^{k-1}\mathcal{H}^{k-1}_{j}(x,t).\]
Recall that by the properties of $M^{k}_{j}$'s:
\begin{itemize}
 \item $M^{k-1}_{j}(T_{[0,t)}(\Phi),C)=M^{k}_{j}(T_{[0,t)}(\Phi),C)$ for $j <\lfloor\frac{k-1}{2}\rfloor$.
 \item $M^{k-1}_{j}(T_{[0,t)}(\Phi),C)=M^{k}_{j+1}(T_{[0,t)}(\Phi),C)$ for $j > \lfloor\frac{k-1}{2}\rfloor$. 
 \item $M^{k-1}_{\lfloor\frac{k-1}{2}\rfloor}(T_{[0,t)}(\Phi),C)=M^{k}_{\lfloor\frac{k-1}{2}\rfloor}(T_{[0,t)}(\Phi),C)\cup M^{k}_{\lfloor\frac{k+1}{2}\rfloor}(T_{[0,t)}(\Phi),C)$.
\end{itemize}
We have then:
\begin{itemize}
 \item $\mathcal{H}^{k-1}_{j}(x,t)=\mathcal{H}^{k}_{j}(x,t)$ for $j< \lfloor\frac{k-1}{2}\rfloor$.
 \item $\mathcal{H}^{k-1}_{j}(x,t)=\mathcal{H}^{k}_{j+1}(x,t)$ for $j> \lfloor\frac{k-1}{2}\rfloor$.
 \item $\mathcal{H}^{k-1}_{\lfloor\frac{k-1}{2}\rfloor}(x,t)=\mathcal{H}^{k}_{\lfloor\frac{k-1}{2}\rfloor}(x,t)\mathcal{H}^{k}_{\lfloor\frac{k+1}{2}\rfloor}(x,t)$.
\end{itemize}
Hence:
\begin{align*}
&\frac{dg_{k,\Lambda}(t,v_{0},\cdots,v_{k})}{dt} =\lim_{\epsilon\rightarrow 0}(g_{k,\Lambda}(t+\epsilon,v_{0},\cdots,v_{k})-g_{k,\Lambda}(t,v_{0},\cdots,v_{k}))\epsilon^{-1}=\lim_{\epsilon\rightarrow 0}\\ 
&\textbf{E}\left[\prod_{i=0}^{k}\prod_{x \in M^{k}_{i}(T_{[0,t)}(\Phi),C)}v_{i}(x)\left(\prod_{i=0}^{k}\prod_{x \in T_{[t,t+\epsilon)}(M^{k}_{i}(T_{[0,t+\epsilon)}(\Phi),C))}v_{i}(x)-1\right)\epsilon^{-1}\right].
\end{align*}
We have:
\begin{align*}
 &\textbf{E}\left[\prod_{i=0}^{k}\prod_{x \in T_{[t,t+\epsilon)}(M^{k}_{i}(T_{[0,t+\epsilon)}(\Phi),C))}v_{i}(x)|T_{[0,t)}(\Phi)\right]=\\
 &\prod_{i=0}^{k}\textbf{E}\left[\prod_{x \in T_{[t,t+\epsilon)}(M^{k}_{i}(T_{[0,t+\epsilon)}(\Phi),C))}v_{i}(x)|T_{[0,t)}(\Phi)\right]=\\
 &\prod_{i=0}^{k-1}\exp\left\{-\epsilon\int_{\mathbb{R}^{2}}(1-v_{i}(x,t))\left(1-\mathcal{H}^{k-1}_{k-1-i}(x)\right)\prod_{j=k-i}^{k-1}\mathcal{H}^{k-1}_{j}(x,t)\Lambda(dx)\right\}\\
 &\exp\left\{-\epsilon\int_{\mathbb{R}^{2}}(1-v_{k}(x))\prod_{j=0}^{k-1}\mathcal{H}^{k-1}_{j}(x,t)\Lambda(dx)\right\}\\
 &= \exp\left\{-\epsilon\left(\sum_{i=0}^{k-1}\int_{\mathbb{R}^{2}}(1-v_{i}(x,t))\left(1-\mathcal{H}^{k-1}_{k-1-i}(x)\right)\prod_{j=k-i}^{k-1}\mathcal{H}^{k-1}_{j}(x,t)\Lambda(dx)\right.\right.\\
 &\left.\left.+\int_{\mathbb{R}^{2}}(1-v_{k}(x))\prod_{j=0}^{k-1}\mathcal{H}^{k-1}_{j}(x,t)\Lambda(dx)\right)\right\}\\
 &= \exp\left\{-\epsilon\left(\sum_{i=0}^{k-1}\int_{\mathbb{R}^{2}}(v_{i}(x)-v_{i+1}(x))\prod_{j=k-1-i}^{k-1}\mathcal{H}^{k-1}_{j}(x,t)\Lambda(dx) +\int_{\mathbb{R}^{2}}(1-v_{0}(x))\Lambda(dx)\right)\right\}\\
 &= \exp\left\{-\epsilon\left(\sum_{i=0}^{\lfloor\frac{k}{2}\rfloor-1}\int_{\mathbb{R}^{2}}(v_{i}(x)-v_{i+1}(x))\prod_{j=k-i}^{k}\mathcal{H}^{k}_{j}(x,t)\Lambda(dx)\right.\right.\\
 &\left.\left.+\sum_{i=\lfloor\frac{k}{2}\rfloor}^{k-1}\int_{\mathbb{R}^{2}}(v_{i}(x)-v_{i+1}(x))\prod_{j=k-i-1}^{k}\mathcal{H}^{k}_{j}(x,t)\Lambda(dx)+\int_{\mathbb{R}^{2}}(1-v_{0}(x))\Lambda(dx)\right)\right\}.
\end{align*}
In the last equality we have used the fact that $\lfloor\frac{k}{2}\rfloor+\lfloor\frac{k-1}{2}\rfloor=k-1$.
Using this we get:
\begin{align*}
 &\frac{dg_{k,\Lambda}(t,v_{0},\cdots,v_{k})}{dt} =\lim_{\epsilon\rightarrow 0}\textbf{E}\left[\prod_{i=0}^{k}\prod_{x \in M^{k}_{i}(T_{[0,t)}(\Phi),C)}v_{i}(x)\left(\exp\left\{-\epsilon\left(\int_{\mathbb{R}^{2}}(1-v_{0}(x))\right.\right.\right.\right.\\
 &\left.\left.\left.\left.\Lambda(dx)+\sum_{i=0}^{\lfloor\frac{k}{2}\rfloor-1}\int_{\mathbb{R}^{2}}(v_{i}(x)-v_{i+1}(x))\prod_{j=k-i}^{k}\mathcal{H}^{k}_{j}(x,t)\Lambda(dx)+\sum_{i=\lfloor\frac{k}{2}\rfloor}^{k-1}\int_{\mathbb{R}^{2}}(v_{i}(x)-v_{i+1}(x))\right.\right.\right.\right.\\
 &\left.\left.\left.\left.\prod_{j=k-i-1}^{k}\mathcal{H}^{k}_{j}(x,t)\Lambda(dx)\right)\right\}-1\right)\epsilon^{-1}\right]\\
 &=- \textbf{E}\left[\prod_{j=0}^{k}\prod_{y \in M^{k}_{j}(T_{[0,t)}(\Phi),C)}v_{j}(y)\left(\int_{\mathbb{R}^{2}}(1-v_{0}(x))\Lambda(dx)+\sum_{i=0}^{\lfloor\frac{k}{2}\rfloor-1}\int_{\mathbb{R}^{2}}(v_{i}(x)-v_{i+1}(x))\right.\right.\\
&\left.\left.\prod_{j=k-i}^{k}\mathcal{H}^{k}_{j}(x,t)\Lambda(dx)+\sum_{i=\lfloor\frac{k}{2}\rfloor}^{k-1}\int_{\mathbb{R}^{2}}(v_{i}(x)-v_{i+1}(x))\prod_{j=k-i-1}^{k}\mathcal{H}^{k}_{j}(x,t)\Lambda(dx)\right)\right]\\
 &= -\int_{\mathbb{R}^{2}}(1-v_{0}(x))\textbf{E}\left[\prod_{j=0}^{k}\prod_{y \in M^{k}_{j}(T_{[0,t)}(\Phi),C)}v_{j}(y)\right]\Lambda(dx)-\sum_{i=0}^{\lfloor\frac{k}{2}\rfloor-1}\int_{\mathbb{R}^{2}}(v_{i}(x)-v_{i+1}(x))\\
&\textbf{E}\left[\prod_{j=0}^{k}\prod_{y \in M^{k}_{j}(T_{[0,t)}(\Phi),C)}v_{j}(y)\prod_{j=k-i}^{k}\mathcal{H}^{k}_{j}(x,t)\right]\Lambda(dx)-\sum_{i=\lfloor\frac{k}{2}\rfloor}^{k-1}\int_{\mathbb{R}^{2}}(v_{i}(x)-v_{i+1}(x))\\
&\textbf{E}\left[\prod_{j=0}^{k}\prod_{y \in M^{k}_{j}(T_{[0,t)}(\Phi),C)}v_{j}(y)\prod_{j=k-i-1}^{k}\mathcal{H}^{k}_{j}(x,t)\right]\Lambda(dx).
\end{align*}
 By noting that
\begin{align*}
&\prod_{y \in M^{k}_{j}(T_{[0,t)}(\Phi),C)}v_{j}(y)\mathcal{H}^{k}_{j}(x,t)=\prod_{y \in M^{k}_{j}(T_{[0,t)}(\Phi),C)}H(v_{j},x)(y),
\end{align*}
we have the wanted result.
\begin{flushright}
 $\Box$
\end{flushright}
\section{ Proof of Proposition \ref{P:Palm gen Matern general}}\label{A:Proof P Palm general}
All the technical conditions to apply Proposition \ref{P:Gen func Palm} have been verified in the proof of Proposition \ref{P:Palm gen infty Matern} and will be omitted here.
We start with the computation of the factorial moment densities. We have, for all $n$-tuple $\textbf{B}^{n}$ of pairwise disjoint Borel sets:
\begin{align*}
& f_{\infty,\Lambda}\left(t,1+\sum_{i=1}^{n}s_{i}\textbf{1}_{B_{i}}\right)\\
&=1+\frac{1}{n!}\sum_{j=1}^{\infty}\int_{(\mathbb{R}^{2})^{j}}\prod_{l=1}^{j}\left(\sum_{i=1}^{n}s_{i}\textbf{1}_{B_{i}}(x_{l})\right)m^{n}_{\infty,\Lambda}(t,\textbf{x}^{n})\Lambda^{i}(d\textbf{x}^{i}).
\end{align*}
Thus:
\begin{align*}
& \frac{d^{n}}{d\textbf{s}^{n}}f_{\infty,\Lambda}(t,1+\sum_{i=1}^{n}s_{i}\textbf{1}_{B_{i}})\Bigr|_{\textbf{s}^{n}=0}\\
&= \frac{d^{n}}{d\textbf{s}^{n}}\int_{(\mathbb{R}^{2})^{n}}\prod_{l=1}^{n}(\sum_{i=1}^{n}s_{i}\textbf{1}_{B_{i}}(x_{l}))m^{n}_{\infty,\Lambda}(t,\textbf{x}^{n})\Lambda^{n}(d\textbf{x}^{n})\Bigr|_{\textbf{s}^{n}=0}.
\end{align*}
We can prove by induction on $n$ that:
\begin{align*}
 \frac{d^{n}}{d\textbf{s}^{n}}\prod_{l=1}^{n}(\sum_{i=1}^{n}s_{i}\textbf{1}_{B_{i}}(x_{l}))=n!\prod_{l=1}^{n}\textbf{1}_{B_{l}}(x_{l}).
\end{align*}
Hence:
\begin{align}
& \frac{d^{n}}{d\textbf{s}^{n}}f_{\infty,\Lambda}(t,1+\sum_{i=1}^{n}s_{i}\textbf{1}_{B_{i}})\Bigr|_{\textbf{s}^{n}=0}= \int_{\textbf{B}^{n}}m^{n}_{\infty,\Lambda}(t,\textbf{x}^{n})\Lambda^{n}(d\textbf{x}^{n}). \label{E:n fac moment}
\end{align}
Now, by Theorem \ref{T:CSMA eq}, we have:
\begin{align*}
& f_{\infty,\Lambda}(t,1+\sum_{i=1}^{n}s_{i}\textbf{1}_{B_{i}})\\
&= 1+\int_{0}^{t}\int_{\mathbb{R}^{2}}f_{\infty,\Lambda}(\tau,H(1,x)-\sum_{i=1}^{n}s_{i}H(\textbf{1}_{B_{i}},x))(\sum_{i=1}^{n}s_{i}\textbf{1}_{B_{i}}(x))\Lambda(dx).
\end{align*}
Subtituting this into (\ref{E:n fac moment}), we get:
\begin{align*}
& \int_{\textbf{B}^{n}}m^{n}_{\infty,\Lambda}(t,\textbf{x}^{n})\Lambda^{n}(d\textbf{x}^{n})\\
&=\frac{d^{n}}{d\textbf{s}^{n}}\int_{0}^{t}\int_{\mathbb{R}^{2}}f_{\infty,\Lambda}\left(\tau,H(1,x)+\sum_{i=1}^{n}s_{i}H(\textbf{1}_{B_{i}},x)\right)\left(\sum_{i=1}^{n}s_{i}\textbf{1}_{B_{i}}(x)\right)\Lambda(dx)\Bigr|_{\textbf{s}^{n}=0}\\
&=\int_{0}^{t}\int_{\mathbb{R}^{2}}\frac{d^{n}}{d\textbf{s}^{n}}f_{\infty,\Lambda}\left(\tau,H(1,x)+\sum_{i=1}^{n}s_{i}H(\textbf{1}_{B_{i}},x)\right)\left(\sum_{i=1}^{n}s_{i}\textbf{1}_{B_{i}}(x)\right)\Bigr|_{\textbf{s}^{n}=0}\Lambda(dx).
\end{align*}
Note that the condition for this interchange of integration and differentiation has beeen established in the proof of Proposition \ref{P:Palm gen infty Matern}.\\
As for all $k>2$, $\frac{d^{k}}{ds_{j_{1}}\cdots ds_{j_{k}}}(\sum_{i=1}^{n}s_{i}\textbf{1}_{B_{i}}(x))\bigr|_{\textbf{s}^{n}=0}=0$ if $\{j_{1},\cdots,j_{k}\}$ is a subsequence of $\{1,\cdots,n\}$, and $(\sum_{i=1}^{n}s_{i}\textbf{1}_{B_{i}})\Bigr|_{\textbf{s}^{n}=0}=0$, we have:
\begin{align*}
 &\frac{d^{n}}{d\textbf{s}^{n}}f_{\infty,\Lambda}\left(\tau,H(1,x)+\sum_{i=1}^{n}s_{i}H(\textbf{1}_{B_{i}},x)\right)\left(\sum_{i=1}^{n}s_{i}\textbf{1}_{B_{i}}(x)\right)\Bigr|_{\textbf{s}^{n}=0}\\
&=\sum_{j=1}^{n}\frac{d^{n-1}}{d\textbf{s}^{n}_{-j}}f_{\infty,\Lambda}\left(\tau,H(1,x)+\sum_{i=1}^{j-1}s_{i}H(\textbf{1}_{B_{i}},x)+\sum_{i=j+1}^{n}s_{i}H(\textbf{1}_{B_{i}},x)\right)\Bigr|_{\textbf{s}^{n}_{-j}=0}\textbf{1}_{B_{j}}(x).
\end{align*}
In the above computation, we used the formula: 
\begin{align}
\frac{d^{n}}{d\textbf{s}^{n}}uv=\sum_{D \subset \llbracket1,n\rrbracket}\frac{d^{|D|}}{\prod_{j \in D}ds_{j}}u\frac{d^{n-|D|}}{\prod_{j \in \llbracket1,n\rrbracket\setminus D}ds_{j}}v.\label{E:multi diff prod}
\end{align}
We use the same argument as with (\ref{E:n fac moment}) to obtain:
\begin{align*}
& \frac{d^{n-1}}{d\textbf{s}^{n}_{-j}}f_{\infty,\Lambda}\left(\tau,H(1,x)\sum_{i \in\llbracket 1,n\rrbracket, i \neq j}s_{i}H(\textbf{1}_{B_{i}},x)\right)\Biggr|_{\textbf{s}^{n}_{-j}=0}\\
&=\frac{d^{n-1}}{d\textbf{s}^{n}_{-j}}f_{\infty,\Lambda}\left(\tau,H(1,x_{j})+\sum_{i \in\llbracket 1,n\rrbracket, i \neq j}s_{i}H(\textbf{1}_{B_{i}},x_{j})\right)\Biggr|_{\textbf{s}^{n}_{-j}=0}\\
&=\int_{\textbf{B}^{n}_{-j}}f^{\textbf{x}^{n}_{-j}}_{\infty,\Lambda}(\tau,H(1,x_{j}))\prod_{i \in\llbracket 1,n\rrbracket, i \neq j}(1-h(x_{i},x_{j}))m^{n-1}_{\infty,\Lambda}(\tau,\textbf{x}^{n}_{-j})\Lambda^{n-1}(d\textbf{x}^{n}_{-j}).
\end{align*}
Hence:
\begin{align*}
&\int_{\textbf{B}^{n}}m^{n}_{\infty,\Lambda}(t,\textbf{x}^{n})\Lambda^{n}(d\textbf{x}^{n})\\
\\
&=\sum_{j=1}^{n}\int_{B_{j}}\int_{\textbf{B}^{n}_{-j}}f^{\textbf{x}^{n}_{-j}}_{\infty,\Lambda}(\tau,H(1,x_{j}))\prod_{i \in\llbracket 1,n\rrbracket, i \neq j}(1-h(x_{i},x_{j}))m^{n-1}_{\infty,\Lambda}(\tau,\textbf{x}^{n}_{-j})\\
&\Lambda^{n-1}(d\textbf{x}^{n}_{-j})\Lambda(dx_{j})\\
&=\int_{\textbf{B}^{n}}f^{\textbf{x}^{n}_{-j}}_{\infty,\Lambda}(\tau,H(1,x_{j}))\prod_{i \in\llbracket 1,n\rrbracket, i \neq j}(1-h(x_{i},x_{j}))m^{n-1}_{\infty,\Lambda}(\tau,\textbf{x}^{n}_{-i})\Lambda^{n}(d\textbf{x}^{n}).
\end{align*}
As this is true for any $n$ tuple $\textbf{B}^{n}$ of disjoint Borel sets, we must have $\Lambda^{n}$-almost everywhere that:
\begin{align*}
 m^{n}_{\infty,\Lambda}(t,\textbf{x}^{n})= \sum_{j=1}^{n}\int_{0}^{t}f^{\textbf{x}^{n}_{-j}}_{\infty,\Lambda}(\tau,H(1,x_{i}))\prod_{i \in\llbracket 1,n\rrbracket,i \neq j}(1-h(x_{i},x_{j}))m^{n-1}_{\infty,\Lambda}(\tau,\textbf{x}^{n}_{-i})d\tau.
\end{align*}
Now, we move to the computation of $f^{\textbf{y}^{n}}_{\infty,\Lambda}(t,v)$. Using the same argument as in (\ref{E:n fac moment}), we have for all $n$-tuple $\textbf{B}^{n}$ of disjoint Borel sets:
\begin{align}
& \frac{d^{n}}{d\textbf{s}^{n}}f_{\infty,\Lambda}(t,v+\sum_{i=1}^{n}s_{i}\textbf{1}_{B_{i}})\Bigr|_{\textbf{s}^{n}=0}= \int_{\textbf{B}^{n}}f^{\textbf{y}^{n}}_{\infty,\Lambda}(t,v)m^{n}_{\infty,\Lambda}(t,\textbf{y}^{n})\Lambda^{n}(d\textbf{y}^{n}). \label{E:n fold Palm}
\end{align}
By using again Theorem \ref{T:CSMA eq}, we get:
\begin{align*}
& f_{\infty,\Lambda}(t,v+\sum_{i=1}^{n}s_{i}\textbf{1}_{B_{i}})= 1-\int_{0}^{t}\int_{\mathbb{R}^{2}}f_{\infty,\Lambda}\left(\tau,H(v,x)-\sum_{i=1}^{n}s_{i}H(\textbf{1}_{B_{i}},x)\right)\\
&\left(1-v(x)-\sum_{i=1}^{n}s_{i}\textbf{1}_{B_{i}}(x)\right)\Lambda(dx)d\tau. 
\end{align*}
 Hence:
\begin{align*}
& \int_{\textbf{B}^{n}}f^{\textbf{y}^{n}}_{\infty,\Lambda}(t,v)m^{n}_{\infty,\Lambda}(t,\textbf{y}^{n})\Lambda^{n}(d\textbf{y}^{n})=-\int_{0}^{t}\int_{\mathbb{R}^{2}}\frac{d^{n}}{d\textbf{s}^{n}}\\
&f_{\infty,\Lambda}\left(\tau,H(v,x)+\sum_{i=1}^{n}s_{i}H(\textbf{1}_{B_{i}},x)\right)\left(1-v(x)-\sum_{i=1}^{n}s_{i}\textbf{1}_{B_{i}}(x)\right)\Biggr|_{\textbf{s}^{n}=0}\Lambda(dx)d\tau.
\end{align*}
By using again (\ref{E:multi diff prod}) we have:
\begin{align*}
 &-\frac{d^{n}}{d\textbf{s}^{n}}f_{\infty,\Lambda}\left(\tau,H(v,x)+\sum_{i=1}^{n}s_{i}H(\textbf{1}_{B_{i}},x)\right)\left(1-v(x)-\sum_{i=1}^{n}s_{i}\textbf{1}_{B_{i}}(x)\right)\Biggr|_{\textbf{s}^{n}=0}\\
 &=\sum_{j=1}^{n}\frac{d^{n-1}}{d\textbf{s}^{n}_{-j}}f_{\infty,\Lambda}\left(\tau,H(v,x)+\sum_{i=1}^{n}s_{i}H(\textbf{1}_{B_{i}},x)\right)\Biggr|_{\textbf{s}^{n}=0}\textbf{1}_{B_{j}}(x)\\
& -\frac{d^{n}}{d\textbf{s}^{n}}f_{\infty,\Lambda}\left(\tau,H(v,x)+\sum_{i=1}^{n}s_{i}H(\textbf{1}_{B_{i}},x)\right)\Biggr|_{\textbf{s}^{n}=0}(1-v(x)).
\end{align*}
By Proposition \ref{P:Gen func Palm}, we have:
\begin{align*}
& \frac{d^{n-1}}{d\textbf{s}^{n}_{-j}}f_{\infty,\Lambda}\left(\tau,H(v,x)+\sum_{i=1}^{n}s_{i}H(\textbf{1}_{B_{i}},x)\right)\Biggr|_{\textbf{s}^{n}=0}\\
&= \frac{d^{n-1}}{d\textbf{s}^{n}_{-j}}f_{\infty,\Lambda}\left(\tau,H(v,y_{j})+\sum_{i=1}^{n}s_{i}H(\textbf{1}_{B_{i}},y_{j})\right)\Biggr|_{\textbf{s}^{n}=0}\\
&=\int_{\textbf{B}^{n}_{-j}}f^{\textbf{y}^{n}_{-j}}_{\infty,\Lambda}(t,v)\prod_{i \in \llbracket1,n\rrbracket, i\neq j}(1-h(y_{i},y_{j}))m^{n-1}_{\infty,\Lambda}(\tau,\textbf{y}^{n}_{-j})\Lambda^{n-1}(d\textbf{y}^{n}_{-j}),\\
&\frac{d^{n}}{d\textbf{s}^{n}}f_{\infty,\Lambda}\left(\tau,H(v,x)+\sum_{i=1}^{n}s_{i}H(\textbf{1}_{B_{i}},x)\right)\Biggr|_{\textbf{s}^{n}=0}\\
&=\int_{\textbf{B}^{n}}f^{\textbf{y}^{n}}_{\infty,\Lambda}(\tau,H(v,x))\prod_{i=1}^{n}(1-h(y_{i},x))(1-v(x))m^{n}_{\infty,\Lambda}\Lambda^{n}(d\textbf{y}^{n}).
\end{align*}
Hence:
\begin{align*}
& \int_{\textbf{B}^{n}}f^{\textbf{y}^{n}}_{\infty,\Lambda}(t,v)m^{n}_{\infty,\Lambda}(t,\textbf{y}^{n})\Lambda^{n}\\
&=\sum_{j=1}^{n}\int_{0}^{t}\int_{\textbf{B}^{n}}f^{\textbf{y}^{n}_{-j}}_{\infty,\Lambda}(t,v)\prod_{i \in \llbracket1,n\rrbracket, i\neq j}(1-h(y_{i},y_{j}))m^{n-1}_{\infty,\Lambda}(\tau,\textbf{y}^{n}_{-j})\Lambda^{n}(d\textbf{y}^{n})d\tau-\\
&\int_{0}^{t}\int_{\textbf{B}^{n}}\int_{\mathbb{R}^{2}}f^{\textbf{y}^{n}}_{\infty,\Lambda}(\tau,H(v,x))\prod_{i=1}^{n}(1-h(y_{i},x))(1-v(x))m^{n}_{\infty,\Lambda}\Lambda^{n}(d\textbf{y}^{n})\Lambda(dx)d\tau. 
\end{align*}
As this is true for all $n$-tuple $\textbf{B}^{n}$ of disjoint Borel sets, we must have $\Lambda^{n}$ almost everywhere:
\begin{align*}
 &f^{\textbf{y}^{n}}_{\infty}(t,v)=\sum_{j=1}^{n}\int_{0}^{t}f^{\textbf{y}^{n}_{-i}}_{\infty,\Lambda}(\tau,H(v,y_{i}))\prod_{i \in\llbracket 1,n\rrbracket,i \neq j}(1-h(y_{i},y_{j}))\frac{m^{n-1}_{\infty,\Lambda}(\tau,\textbf{y}^{n}_{-j})}{m^{n}_{\infty,\Lambda}(t,\textbf{y}^{n})}d\tau\\
&-\int_{0}^{t}\int_{\mathbb{R}^{2}}f^{\textbf{y}^{n}}_{\infty,\Lambda}(\tau,H(v,x))\prod_{i=1}^{n}(1-h(y_{i},x))(1-v(x))\frac{m^{n}_{\infty,\Lambda}(\tau,\textbf{y}^{n})}{m^{n}_{\infty,\Lambda}(t,\textbf{y}^{n})}\Lambda(dx)d\tau.
\end{align*}
\begin{flushright}
 $\Box$.
\end{flushright}
\end{document}